# AVERAGING OF HAMILTONIAN FLOWS WITH AN ERGODIC COMPONENT

By Dmitry Dolgopyat[1] and Leonid Koralov[2]

*University of Maryland*

We consider a process on $\mathbb{T}^2$, which consists of fast motion along the stream lines of an incompressible periodic vector field perturbed by white noise. It gives rise to a process on the graph naturally associated to the structure of the stream lines of the unperturbed flow. It has been shown by Freidlin and Wentzell [*Random Perturbations of Dynamical Systems*, 2nd ed. Springer, New York (1998)] and [*Mem. Amer. Math. Soc.* **109** (1994)] that if the stream function of the flow is periodic, then the corresponding process on the graph weakly converges to a Markov process. We consider the situation where the stream function is not periodic, and the flow (when considered on the torus) has an ergodic component of positive measure. We show that if the rotation number is Diophantine, then the process on the graph still converges to a Markov process, which spends a positive proportion of time in the vertex corresponding to the ergodic component of the flow.

**1. Introduction.** Consider the following stochastic differential equation:

$$(1) \qquad dX_t^\varepsilon = \frac{1}{\varepsilon} v(X_t^\varepsilon)\, dt + dW_t, \qquad X_t^\varepsilon \in \mathbb{T}^2.$$

Here, $v(x)$ is an incompressible periodic vector field, $W_t$ is a two-dimensional Brownian motion and $\varepsilon$ is a small parameter. For simplicity of notation, assume that the period of $v$ in each of the variables is equal to one and that $v$ is infinitely smooth. Let $H(x_1, x_2)$ be the stream function of the flow, that is,

$$\nabla^\perp H = (-H'_{x_2}, H'_{x_1}) = v.$$

Received February 2006; revised May 2007.

[1]Supported in part by the NSF and IPST.

[2]Supported in part by the NSF.

Since $v$ is periodic, we can write $H$ as

$$H(x_1, x_2) = H_0(x_1, x_2) + ax_1 + bx_2,$$

where $H_0$ is periodic. We shall assume that all the critical points of $H$ are nondegenerate, and that $(a, b)$ satisfy the following Diophantine condition.

Let $\rho = a/b$ be irrational. Without loss of generality, we may assume that $0 < a < b$ (the general case can be obtained by interchanging $x_1$ and $x_2$, and/or replacing $x_i$ by $-x_i$, if needed). Let $[a_1, a_2 \ldots a_n \ldots]$ be the continued fraction expansion of $\rho$, We assume that

$$(2) \qquad\qquad a_n \leq n^2 \qquad \text{for all sufficiently large } n.$$

We shall see in Section A.3 that this condition holds for almost all $\rho$ (with respect to the Lebesgue measure on $[0, 1]$).

For $a$ and $b$ which are rationally independent, as in our case, it has been shown by Arnold in [1] that the structure of the streamlines of $v$ on the torus is as follows. There are finitely many domains $U_k$, $k = 1, \ldots, n$, bounded by the separatrices of $H$, such that the trajectories of the dynamical system $\dot{X}_t = v(X_t)$ in $U_k$ are either periodic or tend to a point where the vector field is equal to zero. The trajectories form one ergodic class outside the domains $U_k$. More precisely, let $\mathcal{E} = \mathbb{T}^2 \setminus \mathrm{Cl}(\bigcup_{k=1}^n U_k)$. Here, $\mathrm{Cl}(\cdot)$ stands for the closure of a set. The dynamical system is then ergodic on $\mathcal{E}$ (and is, in fact, mixing; see [8]).

Although $H$ itself is not periodic, we can consider its critical points as points on the torus since $\nabla H$ is periodic. All the maxima and minima of $H$ are located inside the domains $U_k$.

At first, we shall consider the case where there is just one periodic component $U$, which contains only one critical point of $H$ (a maximum or a minimum). An example of a phase portrait of such a vector field $v$ (considered on the plane) is given in Figure 1. The general case is discussed in Section A.4.

Assume, for definiteness, that the critical point of $H$ inside $U$ is a maximum. We shall denote the saddle point of $H$ on the torus by $A$ and the maximum by $M$. Consider the following mapping of the torus onto the segment $I = [0, H(M) - H(A)]$ of the real line:

$$h(x) = \begin{cases} 0, & \text{if } x \in \mathcal{E}, \\ H(x) - H(A), & \text{otherwise.} \end{cases}$$

We denote the set $\{x \in \mathrm{Cl}(U) : H(x) - H(A) = h\}$ by $\gamma(h)$. Let $\gamma = \gamma(0) = \partial U$. Let $Lf = a(h)f'' + b(h)f'$ be the differential operator on $I$ with the coefficients

$$(3) \qquad\qquad a(h) = \frac{1}{2} \left( \int_{\gamma(h)} \frac{1}{|\nabla H|} \, dl \right)^{-1} \int_{\gamma(h)} |\nabla H| \, dl$$

and

$$b(h) = \frac{1}{2} \left( \int_{\gamma(h)} \frac{1}{|\nabla H|} \, dl \right)^{-1} \int_{\gamma(h)} \frac{\Delta H}{|\nabla H|} \, dl. \tag{4}$$

Let $k = 2(\int_\gamma |\nabla H| \, dl)^{-1} \operatorname{Area}(\mathcal{E})$. Consider the process $Y_t$ on the segment $I$ which is defined via its generator $\mathcal{L}$ as follows. The domain of the generator $D(\mathcal{L})$ consists of those functions $f \in C(I)$ which:

    (a) are twice continuously differentiable in the interior of $I$;

    (b) have limits $\lim_{h \to 0} Lf(h)$ and $\lim_{h \to (H(M)-H(A))} Lf(h)$ at the endpoints of $I$;

    (c) have the limit $\lim_{h \to 0} f'(h)$ and $\lim_{h \to 0} f'(h) = k \lim_{h \to 0} Lf(h)$.

For functions $f$ which satisfy the above three properties, we define $\mathcal{L}f = Lf$ in the interior of the segment and as the limit of $Lf$ at the endpoints of $I$.

It is well known (see, e.g., [11]) that there exists a strong Markov process on $I$ with continuous trajectories, with the generator $\mathcal{L}$. The measure on $C([0,\infty), I)$ induced by the process is uniquely defined by the operator and the initial distribution of the process.

We shall prove the following theorem.

THEOREM 1. *The measure on $C([0,\infty), I)$ induced by the process $Y_t^\varepsilon = h(X_t^\varepsilon)$ converges weakly to the measure induced by the process with the generator $\mathcal{L}$ with the initial distribution $h(X_0^\varepsilon)$.*

One of the main ingredients of the proof is the estimate of the expectation of the time it takes for the solution of (1) to leave the ergodic component. This estimate will be derived in Sections 4 and 5. Besides this, we shall use a number of estimates on the transition times of the process between different level sets of the Hamiltonian inside the periodic component. Those will be proven in Section 3. First, however, we prove Theorem 1, while assuming that we have all the needed estimates. We discuss the case of more than one periodic component in Section A.4. Sections A.1–A.3 contain technical estimates needed for the proof.

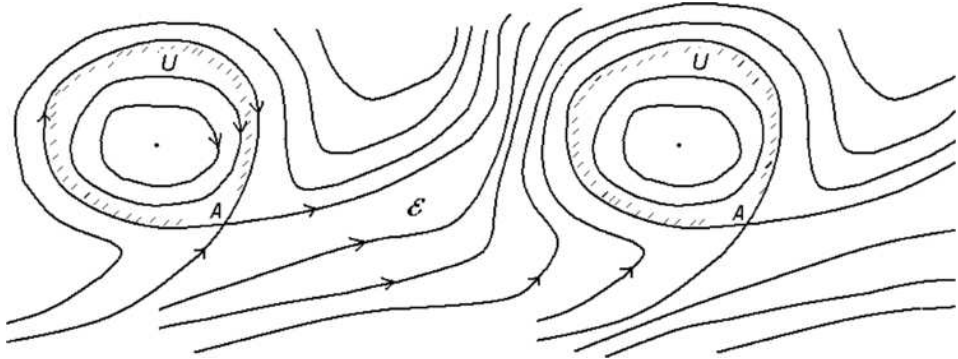

FIG. 1. *Structure of the stream lines.*

We observe that irrationality of $\rho$ is necessary for Theorem 1 since if $\rho$ is rational, then the restriction of $X$ to $\mathcal{E}$ is periodic rather than ergodic, so the phase space of the limiting process is a graph with two edges (one of which forms a loop) and one vertex (cf. [7]). On the other hand, it has been conjectured by Freidlin [5] that Theorem 1 holds for *all* irrational values of $\rho$. Our paper proves this result for $\rho$'s which cannot be approximated too well by rationals. On the other hand, Sowers [14] shows that in the opposite case of $\rho$'s which are very well approximable by rationals, the result is also true. There is still a gap between our condition (2) and the numbers considered in [14], but we hope that the combination of our approach with the methods of [14] will allow the result to be established in full generality.

We note that there is a related problem, which concerns the asymptotics of effective diffusivity for a two-dimensional periodic vector field perturbed by small diffusion. It also involves the study of the behavior of the process $X_t^{\varepsilon}$ near the saddle points and near the separatrices. We refer the interested reader to the papers by Fannjiang and Papanicolaou [4], Koralov [10], Sowers [15] and Novikov, Papanicolaou and Ryzhik [13] for some of the recent results.

**2. Proof of the main theorem.** Let $\Psi$ be the subset of $C(I)$ which consists of all bounded functions which are continuously differentiable on $[0, H(M) - H(A))$ (the derivative at $h = 0$ is one sided). Note that this is a measure-defining set, that is, the equality $\int_I u \, d\mu_1 = \int_I u \, d\mu_2$ for all $u \in \Psi$ implies that $\mu_1 = \mu_2$. Let $\mathcal{D}$ be the subset of $D(\mathcal{L})$ which consists of all the functions $f$ for which $\mathcal{L}f \in \Psi$.

We formulate the following lemma.

LEMMA 2.1. *For any function $f \in \mathcal{D}$, any initial point $x \in \mathbb{T}^2$ and any $T > 0$, we have*

$$(5) \quad \mathbb{E}_x \left[ f(h(X_T^{\varepsilon})) - f(h(X_0^{\varepsilon})) - \int_0^T \mathcal{L}f(h(X_s^{\varepsilon})) \, ds \right] \to 0 \qquad \text{as } \varepsilon \to 0,$$

*uniformly in $x \in \mathbb{T}^2$.*

An analogous lemma was used in the monograph of Freidlin and Wentzell [7] to justify the convergence of the process $Y_t^{\varepsilon}$ to the limiting process on the graph. The main idea, roughly speaking, is to use the tightness of the family $Y_t^{\varepsilon}$ and then to show that the limiting process (along any subsequence) is a solution of the martingale problem, corresponding to the operator $\mathcal{L}$. Here, as in [7], the fact that for every $u \in \Psi$ and $\lambda > 0$, the equation $\lambda f - \mathcal{L}f = u$ has a solution $f \in \mathcal{D}$ is used.

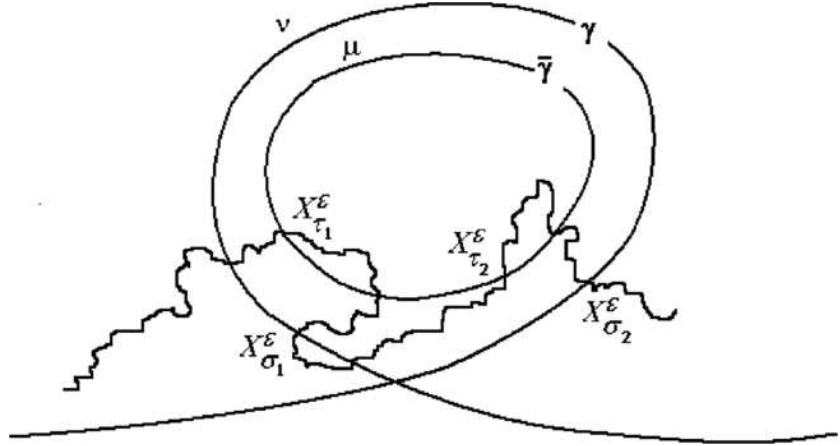

Fig. 2. *Definition of the stopping times.*

The main difference between our case and that of [7] is the presence of an ergodic component. However, all the arguments used to prove the main theorem based on (5) remain the same. Thus, in referring to Lemma 3.1 of [7], it is enough to prove our Lemma 2.1 above.

The proof of Lemma 2.1 will rely on several other lemmas. Below, we shall introduce a number of processes, stopping times and sets, which will depend on $\varepsilon$. However, we shall not always incorporate this dependence on $\varepsilon$ into notation, so one must be careful to distinguish between the objects which do not depend on $\varepsilon$ and those which do.

Let $\overline{\tau}$ be the first time when the process $X_t^{\varepsilon}$ reaches the set $\gamma(\varepsilon^{1/2})$. We will need the following estimate on the expectation of $\overline{\tau}$, which is proved in Sections 4 and 5.

LEMMA 2.2. *For any $\varkappa > 0$, there exists some $\varepsilon_0 > 0$ such that $\mathbb{E}_x \overline{\tau} \leq \varepsilon^{1/2-\varkappa}$ for $\varepsilon \leq \varepsilon_0$, for all $x \in \mathrm{Cl}(\mathcal{E})$.*

Let us now choose constants $\varkappa$ and $\alpha$, such that $0 < \varkappa < \frac{1}{4} < \alpha < \frac{1}{2}$. Let $\overline{\gamma} = \gamma(\varepsilon^{\alpha})$. Recall that $\gamma = \gamma(0)$ is the boundary of $U$. Let $\tau$ be the first time when the process $X_t^{\varepsilon}$ reaches $\overline{\gamma}$ and $\sigma$ be the first time when the process reaches $\gamma$. We inductively define the following two sequences of stopping times. Let $\tau_1 = \tau$. For $n \geq 1$, let $\sigma_n$ be the first time following $\tau_n$ when the process reaches $\gamma$. For $n \geq 2$, let $\tau_n$ be the first time following $\sigma_{n-1}$ when the process reaches $\overline{\gamma}$. See Figure 2.

We consider the discrete-time Markov chains $\xi_n^1 = X_{\tau_n}^\varepsilon$ and $\xi_n^2 = X_{\sigma_n}^\varepsilon$ with the state spaces $\overline{\gamma}$ and $\gamma$, respectively. Let $P_1(x, dy)$ and $P_2(x, dy)$ be transition operators for the Markov chains $\xi_n^1$ and $\xi_n^2$, respectively. They are uniformly exponentially mixing in the following sense.

LEMMA 2.3.   *There exist constants $0 < c < 1$, $\varepsilon_0 > 0$, $n_0 > 0$ and probability measures $\mu$ and $\nu$ (which depend on $\varepsilon$) on $\overline{\gamma}$ and $\gamma$, respectively, such that for $\varepsilon < \varepsilon_0$ and $n \geq n_0$, we have*

$$
\sup_{x \in \overline{\gamma}}(\operatorname{Var}(P_1^n(x, dy) - \mu(dy))) \leq c^n,
$$

(6)

$$
\sup_{x \in \gamma}(\operatorname{Var}(P_2^n(x, dy) - \nu(dy))) \leq c^n,
$$

*where* Var *is the total variation of the signed measure.*

We prove this lemma in Section A.2. Let us now examine the transition times between $\overline{\gamma}$ and $\gamma$, assuming that we start with the invariant measures. The following lemma is proved in Section 3.

LEMMA 2.4.   *The asymptotic behavior of the transition times is the following:*

(7)          $$\mathbb{E}_\mu \sigma = k_1 \varepsilon^\alpha (1 + o(1)) \qquad \text{as } \varepsilon \to 0,$$

(8)          $$\mathbb{E}_\nu \tau = k_2 \varepsilon^\alpha (1 + o(1)) \qquad \text{as } \varepsilon \to 0,$$

*where $k_1 = 2(\int_\gamma |\nabla H| \, dl)^{-1} \operatorname{Area}(U)$ and $k_2 = 2(\int_\gamma |\nabla H| \, dl)^{-1} \operatorname{Area}(\mathcal{E})$. Further, we have the following estimate:*

$$
\sup_{x \in \gamma} \mathbb{E}_x \tau \leq k_3 \varepsilon^{\alpha - \varkappa} \qquad \text{as } \varepsilon \to 0
$$

*for some $k_3 > 0$.*

We will need to control the number of excursions between $\overline{\gamma}$ and $\gamma$ before time $T$. For this purpose, we formulate the following lemma, which will be proven in Section 3.

LEMMA 2.5.   *There is a constant $r > 0$ such that for all sufficiently small $\varepsilon$, we have*

$$
\sup_{x \in \overline{\gamma}} \mathbb{E}_x e^{-\sigma} \leq 1 - r \varepsilon^\alpha.
$$

Using the Markov property of the process and Lemma 2.5, we get the estimate

$$\sup_{x \in \mathbb{T}^2} \mathbb{E}_x e^{-\sigma_n} = \sup_{x \in \overline{\gamma}} \mathbb{E}_x e^{-\sigma_n} \leq \left( \sup_{x \in \overline{\gamma}} \mathbb{E}_x e^{-\sigma} \right)^n \leq (1 - r\varepsilon^\alpha)^n. \tag{9}$$

The next lemma (also proved in Section 3) allows us to estimate integrals of the type (5) over intervals $[0, \tau]$ and $[0, \sigma]$.

LEMMA 2.6. *For any function $f \in \mathcal{D}$, we have the following asymptotic estimates:*

$$\sup_{x \in \overline{\gamma}} \left| \mathbb{E}_x \left[ f(h(X_\sigma^\varepsilon)) - f(h(X_0^\varepsilon)) - \int_0^\sigma \mathcal{L}f(h(X_s^\varepsilon)) \, ds \right] \right| = o(\varepsilon^\alpha) \tag{10}$$

$$as \; \varepsilon \to 0,$$

$$\mathbb{E}_\nu \left[ f(h(X_\tau^\varepsilon)) - f(h(X_0^\varepsilon)) - \int_0^\tau \mathcal{L}f(h(X_s^\varepsilon)) \, ds \right] = o(\varepsilon^\alpha) \tag{11}$$

$$as \; \varepsilon \to 0.$$

*Moreover, we also have*

$$\sup_{x \in \mathbb{T}^2} \left| \mathbb{E}_x \left[ f(h(X_\tau^\varepsilon)) - f(h(X_0^\varepsilon)) - \int_0^\tau \mathcal{L}f(h(X_s^\varepsilon)) \, ds \right] \right| \to 0 \tag{12}$$

$$as \; \varepsilon \to 0.$$

PROOF OF LEMMA 2.1. Let $f \in \mathcal{D}$, $T > 0$ and $\eta > 0$ be fixed. We would like to show that the absolute value of the left-hand side of (5) is less than $\eta$ for all sufficiently small positive $\varepsilon$.

First, we replace the time interval $[0, T]$ by a larger one, $[0, \widetilde{\tau}]$, where $\widetilde{\tau}$ is the first of the stopping times $\tau_n$ which is greater than or equal to $T$, that is,

$$\widetilde{\tau} = \min_{n : \tau_n \geq T} \tau_n.$$

Using the Markov property of the process, the difference can be rewritten as follows:

$$\left| \mathbb{E}_x \left[ f(h(X_{\widetilde{\tau}}^\varepsilon)) - f(h(X_0^\varepsilon)) - \int_0^{\widetilde{\tau}} \mathcal{L}f(h(X_s^\varepsilon)) \, ds \right] \right.$$

$$- \mathbb{E}_x \left[ f(h(X_T^\varepsilon)) - f(h(X_0^\varepsilon)) - \int_0^T \mathcal{L}f(h(X_s^\varepsilon)) \, ds \right] \Bigg|$$

$$= \left| \mathbb{E}_x \mathbb{E}_{X_T^\varepsilon} \left[ f(h(X_\tau^\varepsilon)) - f(h(X_0^\varepsilon)) - \int_0^\tau \mathcal{L}f(h(X_s^\varepsilon)) \, ds \right] \right|.$$

The latter expression can be made smaller than $\frac{\eta}{5}$ for all sufficiently small $\varepsilon$ due to (12). Therefore, it remains to show that

$$\left| \mathbb{E}_x\left[ f(h(X_{\widetilde{\tau}}^{\varepsilon})) - f(h(X_0^{\varepsilon})) - \int_0^{\widetilde{\tau}} \mathcal{L}f(h(X_s^{\varepsilon}))\, ds \right] \right| < \frac{4\eta}{5}$$

for all sufficiently small $\varepsilon$. We shall denote the indicator function of a set $A$ by $\chi_A$. Using the stopping times $\tau_n$ and $\sigma_n$, we can rewrite the expectation in the left-hand side of this inequality as follows:

$$\mathbb{E}_x\left[ f(h(X_{\widetilde{\tau}}^{\varepsilon})) - f(h(X_0^{\varepsilon})) - \int_0^{\widetilde{\tau}} \mathcal{L}f(h(X_s^{\varepsilon}))\, ds \right]$$

$$= \mathbb{E}_x\left[ f(h(X_{\tau}^{\varepsilon})) - f(h(X_0^{\varepsilon})) - \int_0^{\tau} \mathcal{L}f(h(X_s^{\varepsilon}))\, ds \right]$$

$$+ \sum_{n=1}^{\infty} \mathbb{E}_x\left( \chi_{\{\tau_n < \widetilde{\tau}\}} \mathbb{E}_{X_{\tau_n}^{\varepsilon}}\left[ f(h(X_{\sigma}^{\varepsilon})) - f(h(X_0^{\varepsilon})) - \int_0^{\sigma} \mathcal{L}f(h(X_s^{\varepsilon}))\, ds \right] \right)$$

$$+ \sum_{n=1}^{\infty} \mathbb{E}_x\left( \chi_{\{\sigma_n < \widetilde{\tau}\}} \mathbb{E}_{X_{\sigma_n}^{\varepsilon}}\left[ f(h(X_{\tau}^{\varepsilon})) - f(h(X_0^{\varepsilon})) - \int_0^{\tau} \mathcal{L}f(h(X_s^{\varepsilon}))\, ds \right] \right),$$

provided that the sums in the right-hand side converge absolutely (which follows from the arguments below). Due to (12), the absolute value of the first term on the right-hand side of this equality can be made smaller than $\frac{\eta}{5}$ for all sufficiently small $\varepsilon$. Therefore, it remains to estimate the two infinite sums.

Let us start with the first sum. Note that

$$\mathbb{E}_x \chi_{\{\tau_n < \widetilde{\tau}\}} \leq \mathbb{E}_x \chi_{\{\sigma_{n-1} < T\}} \leq \mathbb{E}_x \chi_{\{e^{-\sigma_{n-1}} > e^{-T}\}} \leq e^T (1 - r\varepsilon^{\alpha})^{n-1}.$$

The last inequality here is due to (9) and Chebyshev's inequality. Taking the sum in $n$, we obtain

$$\sum_{n=1}^{\infty} \mathbb{E}_x \chi_{\{\tau_n < \widetilde{\tau}\}} \leq \sum_{n=1}^{\infty} e^T (1 - r\varepsilon^{\alpha})^{n-1} \leq K\varepsilon^{-\alpha},$$

where the constant $K$ depends on $T$ and $r$. By Lemma 2.6, we can find $\varepsilon_0$ such that for all $\varepsilon < \varepsilon_0$, we have

$$\sup_{x \in \overline{\gamma}} \left| \mathbb{E}_x\left[ f(h(X_{\sigma}^{\varepsilon})) - f(h(X_0^{\varepsilon})) - \int_0^{\sigma} \mathcal{L}f(h(X_s^{\varepsilon}))\, ds \right] \right| \leq \frac{\eta \varepsilon^{\alpha}}{5K}.$$

Therefore, for $\varepsilon < \varepsilon_0$, we have

$$\left| \sum_{n=1}^{\infty} \mathbb{E}_x\left( \chi_{\{\tau_n < \widetilde{\tau}\}} \mathbb{E}_{X_{\tau_n}^{\varepsilon}}\left[ f(h(X_{\sigma}^{\varepsilon})) - f(h(X_0^{\varepsilon})) - \int_0^{\sigma} \mathcal{L}f(h(X_s^{\varepsilon}))\, ds \right] \right) \right|$$

$$\leq \sup_{x \in \overline{\gamma}} \left| \mathbb{E}_x \left[ f(h(X_\sigma^\varepsilon)) - f(h(X_0^\varepsilon)) - \int_0^\sigma \mathcal{L} f(h(X_s^\varepsilon)) \, ds \right] \right| \sum_{n=1}^\infty \mathbb{E}_x \chi_{\{\tau_n < \widetilde{\tau}\}}$$

$$\leq \frac{\eta}{5}.$$

The same argument allows us to write the estimate as

$$\left| \sum_{n=1}^\infty \mathbb{E}_x \left( \chi_{\{\sigma_n < \widetilde{\tau}\}} \mathbb{E}_\nu \left[ f(h(X_\tau^\varepsilon)) - f(h(X_0^\varepsilon)) - \int_0^\tau \mathcal{L} f(h(X_s^\varepsilon)) \, ds \right] \right) \right| \leq \frac{\eta}{5}.$$

The left-hand side of this inequality, however, is slightly different from the desired expression, since (11) only allows us to estimate the expectation with respect to the original distribution $\nu$, rather than the supremum over all possible initial points. Thus, we need to estimate the difference

$$\left| \sum_{n=1}^\infty \mathbb{E}_x \left( \chi_{\{\sigma_n < \widetilde{\tau}\}} \mathbb{E}_{X_{\sigma_n}^\varepsilon} \left[ f(h(X_\tau^\varepsilon)) - f(h(X_0^\varepsilon)) - \int_0^\tau \mathcal{L} f(h(X_s^\varepsilon)) \, ds \right] \right) \right.$$

$$\left. - \sum_{n=1}^\infty \mathbb{E}_x \left( \chi_{\{\sigma_n < \widetilde{\tau}\}} \mathbb{E}_\nu \left[ f(h(X_\tau^\varepsilon)) - f(h(X_0^\varepsilon)) - \int_0^\tau \mathcal{L} f(h(X_s^\varepsilon)) \, ds \right] \right) \right|.$$

This expression can be estimated from above by

$$\sup_{x \in \gamma} \left| \mathbb{E}_x \left[ f(h(X_\tau^\varepsilon)) - f(h(X_0^\varepsilon)) - \int_0^\tau \mathcal{L} f(h(X_s^\varepsilon)) \, ds \right] \right|$$

$$\times \sum_{n=1}^\infty \sup_{x \in \gamma} (\text{Var}(P_2^{n-1}(x, dy) - \nu(dy))),$$

which is smaller than $\frac{\eta}{5}$ for all sufficiently small $\varepsilon$, due to (6) and (12). Combining the above estimates, we see that the absolute value of the left-hand side of (5) is less than $\eta$ for all sufficiently small positive $\varepsilon$. This completes the proof of the theorem. $\square$

## 3. Behavior inside the periodic component.

In this section, we prove Lemmas 2.4, 2.5 and 2.6. At several points, we shall use Lemma 2.2, which is proved in Section 4 and 5. This involves no circular reasoning since the results of this section are not used in Sections 4 and 5.

We shall need a number of statements from [7] and [10] which describe the limiting behavior of the process $X_t^\varepsilon$ with the initial point $x \in U$, both in the case when $x$ is fixed, and when $x$ is asymptotically close (as a power of $\varepsilon$) to the boundary of $U$.

It has been shown in [7] (Theorem 2.2) that if $X_t^\varepsilon = x \in U$ is fixed, then the process $h(X_t^\varepsilon)$, stopped at $t = \sigma$, converges weakly to the diffusion process on $I$ with the generator $L$, which starts at $h(x)$ and is stopped at the moment when it reaches zero (the left endpoint of $I$). We formulate the following statements, which easily follow from the proof of Theorem 2.2 of [7], as a separate lemma.

LEMMA 3.1. (a) *For any function $f \in \mathcal{D}$, we have*

$$\lim_{\varepsilon \to 0} \sup_{x \in U} \left| \mathbb{E}_x \left[ f(h(X_\sigma^\varepsilon)) - f(h(X_0^\varepsilon)) - \int_0^\sigma \mathcal{L} f(h(X_s^\varepsilon)) \, ds \right] \right| = 0.$$

(b) *For any $\widetilde{h} > 0$, there is a constant $c(\widetilde{h}) > 0$ such that*

$$\lim_{\varepsilon \to 0} \sup_{x \in \gamma(\widetilde{h})} \mathbb{E}_x e^{-\sigma} \le 1 - c(\widetilde{h}).$$

We will also need the following lemma, which gives us the asymptotics of the time the process needs in order to exit the periodic component if the original point is asymptotically close to the boundary. It was proven in [10] (Lemma 4.4).

LEMMA 3.2. *There is a constant $k_1$ such that for any $\frac{1}{4} < \alpha < \frac{1}{2}$, we have*

$$\lim_{\varepsilon \to 0} \sup_{x \in \gamma(\varepsilon^\alpha)} \left| \frac{\mathbb{E}_x \sigma}{\varepsilon^\alpha} - k_1 \right| = 1.$$

PROOF OF LEMMA 2.4. The first part of the lemma follows from Lemma 3.2. Namely, it states that there exists a constant $k_1$ such that

$$\lim_{\varepsilon \to 0} \sup_{x \in \overline{\gamma}} \left| \frac{\mathbb{E}_x \sigma}{\varepsilon^\alpha} - k_1 \right| = 1. \tag{13}$$

Let us recall how to identify the constant $k_1$ (rigorous arguments can be found in [10]). If $x \in U$ did not depend on $\varepsilon$, then the asymptotics of $\mathbb{E}_x \sigma$ could be obtained using the results of [7]. Namely, recall the definition of the differential operator $L$ from Section 1 and let $u(h)$ be the bounded solution of the ordinary differential equation

$$Lu = a(h)u'' + b(h)u' = -1, \qquad h \in \mathrm{Int}(I), \tag{14}$$

with the boundary condition $u(0) = 0$. Such a solution exists and is unique (see, e.g., [7]). It is equal to the expectation of the time it takes for the limiting process, starting at $h$, to reach the endpoint of $I$ corresponding to the boundary of the periodic component. It was demonstrated in [7] (Lemma 2.3) that $\lim_{\varepsilon \to 0} \mathbb{E}_x \sigma = u(h(x))$. In particular,

$$\lim_{\varepsilon \to 0} \mathbb{E}_x \sigma = u'(0)h(x) + o(h(x)) \qquad \text{as } h(x) \to 0.$$

Formula $(13)$ is the corresponding asymptotic formula in the case where $h(x)$ is a function of $\varepsilon$, that is, $h(x) = \varepsilon^{\alpha}$. In particular, $k_1 = u'(0)$. Equation $(14)$ can be solved explicitly using the expressions for the coefficients $a(h)$ and $b(h)$. We obtain that

$$u'(0) = 2\left(\int_{\gamma} |\nabla H|\, dl\right)^{-1} \int_{I}\left(\int_{\gamma(h)} \frac{1}{|\nabla H|}\, dl\right) dh = 2\left(\int_{\gamma} |\nabla H|\, dl\right)^{-1} \mathrm{Area}(U),$$

which proves $(7)$.

In order to study the asymptotics of $\mathbb{E}_{\nu}\tau$, we shall use the asymptotics of $\mathbb{E}_{\mu}\sigma$ and the fact that $X_t^{\varepsilon}$ is an ergodic process on $\mathbb{T}^2$. Let $V^{\varepsilon} = \{x \in U : 0 \le h(x) \le \varepsilon^{\alpha}\}$.

The process $X_t^{\varepsilon}$ is ergodic and the invariant measure is the Lebesgue measure on $\mathbb{T}^2$. Applying the Birkhoff ergodic theorem to the process with the initial distribution $\nu$, we obtain that

$$\lim_{n \to \infty} \frac{\int_0^n \chi_{\mathcal{E}}(X_t^{\varepsilon})\, dt}{n} = \mathrm{Area}(\mathcal{E}) \qquad \text{almost surely},$$

where $\chi_{\mathcal{E}}$ is the indicator function of the set $\mathcal{E}$. Also, from the Birkhoff ergodic theorem, we get

$$\lim_{n \to \infty} \frac{\sigma_n}{n} = \mathbb{E}_{\nu}\sigma_1 = \mathbb{E}_{\nu}\tau + \mathbb{E}_{\mu}\sigma \qquad \text{almost surely}.$$

Using the Birkhoff ergodic theorem again, we can write

$$\lim_{n \to \infty} \frac{\int_0^n \chi_{\mathcal{E}}(X_t^{\varepsilon})\, dt}{n} = \lim_{n \to \infty} \frac{\int_0^{\sigma_n} \chi_{\mathcal{E}}(X_t^{\varepsilon})\, dt}{\sigma_n} = \left(\lim_{n \to \infty} \frac{\int_0^{\sigma_n} \chi_{\mathcal{E}}(X_t^{\varepsilon})\, dt}{n}\right)\left(\lim_{n \to \infty} \frac{n}{\sigma_n}\right)$$

$$= \frac{\mathbb{E}_{\nu}\int_0^{\sigma_1} \chi_{\mathcal{E}}(X_t^{\varepsilon})\, dt}{\mathbb{E}_{\nu}\tau + \mathbb{E}_{\mu}\sigma} \le \frac{\mathbb{E}_{\nu}\tau}{\mathbb{E}_{\nu}\tau + \mathbb{E}_{\mu}\sigma},$$

where the equalities hold almost surely. Therefore,

$$\mathrm{Area}(\mathcal{E}) \le \frac{\mathbb{E}_{\nu}\tau}{\mathbb{E}_{\nu}\tau + \mathbb{E}_{\mu}\sigma}.$$

In exactly the same way, we can prove that

$$\mathrm{Area}(\mathcal{E} \cup V^{\varepsilon}) \ge \frac{\mathbb{E}_{\nu}\tau}{\mathbb{E}_{\nu}\tau + \mathbb{E}_{\mu}\sigma}, \qquad \mathrm{Area}(U \setminus V^{\varepsilon}) \le \frac{\mathbb{E}_{\mu}\sigma}{\mathbb{E}_{\nu}\tau + \mathbb{E}_{\mu}\sigma}$$

and

$$\mathrm{Area}(U) \ge \frac{\mathbb{E}_{\mu}\sigma}{\mathbb{E}_{\nu}\tau + \mathbb{E}_{\mu}\sigma}.$$

Combining these four estimates, we obtain that

$$\frac{\mathrm{Area}(\mathcal{E})}{\mathrm{Area}(\mathcal{U})} \le \frac{\mathbb{E}_{\nu}\tau}{\mathbb{E}_{\mu}\sigma} \le \frac{\mathrm{Area}(\mathcal{E} \cup V^{\varepsilon})}{\mathrm{Area}(\mathcal{U} \setminus V^{\varepsilon})}.$$

Since $\lim_{\varepsilon \to 0} \mathrm{Area}(V^\varepsilon) = 0$, we obtain

$$\lim_{\varepsilon \to 0} \frac{\mathbb{E}_\nu \tau}{\mathbb{E}_\mu \sigma} = \frac{\mathrm{Area}(\mathcal{E})}{\mathrm{Area}(\mathcal{U})}. \tag{15}$$

Therefore, (7) implies (8). It remains to prove the last statement of the lemma.

Let $\gamma_1 = \{x : h(x) = \varepsilon^{1/2}\}$. We inductively define the following two sequences of stopping times. Let $\overline{\tau}_1 = \overline{\tau}$ be the first time when the process $X_t^\varepsilon$ reaches the set $\gamma_1$. For $n \geq 1$, let $\overline{\sigma}_n$ be the first time following $\overline{\tau}_n$ when the process reaches $\gamma$. For $n \geq 2$, let $\overline{\tau}_n$ be the first time following $\overline{\sigma}_{n-1}$ when the process reaches $\gamma_1$.

First, let us estimate the probability of the event that the process which starts at $x \in \gamma_1$ reaches $\gamma_1$ before reaching $\gamma$. Lemma 4.3 of [10] states that there is a constant $c_1$ such that for any $x \in \overline{\gamma}$,

$$|\mathbb{P}_x(\tau < \sigma) - \varepsilon^{1/2-\alpha}| \leq c_1 \varepsilon^\alpha |\ln \varepsilon|.$$

Since $\alpha > \frac{1}{4}$, this implies that $\mathbb{P}_x(\tau > \sigma) \leq 1 - \frac{1}{2}\varepsilon^{1/2-\alpha}$ for all sufficiently small $\varepsilon$, for all $x \in \gamma_1$. Using the Markov property of the process, we conclude that

$$\sup_{x \in \gamma} \mathbb{P}_x(\tau > \overline{\sigma}_n) \leq (1 - \tfrac{1}{2}\varepsilon^{1/2-\alpha})^n.$$

We also need to estimate how much time it takes for the process which starts at $x \in \gamma_1$ to leave $V^\varepsilon$ (the region between $\gamma$ and $\overline{\gamma}$). Lemma 4.2 of [10] states that there is a constant $c_2$ such that for any $x \in V^\varepsilon$,

$$\mathbb{E}_x \min(\tau, \sigma) \leq c_2 \varepsilon^{2\alpha} |\ln \varepsilon|.$$

Since $2\alpha > 1/2 - \varkappa$, the right-hand side of this inequality is smaller than $\varepsilon^{1/2-\varkappa}$ for all sufficiently small $\varepsilon$. Therefore, by Lemma 2.2,

$$\sup_{x \in \gamma} \mathbb{E}_x \min(\tau, \overline{\sigma}_1) \leq \sup_{x \in \gamma} \mathbb{E}_x \overline{\tau} + \sup_{x \in \overline{\gamma}} \mathbb{E}_x \min(\tau, \sigma) \leq 2\varepsilon^{1/2-\varkappa}.$$

Finally, due to the Markov property of the process,

$$\sup_{x \in \gamma} \mathbb{E}_x \tau = \sup_{x \in \gamma} \left[ \mathbb{E}_x \min(\tau, \overline{\sigma}_1) + \sum_{n=1}^\infty (\mathbb{E}_x \min(\tau, \overline{\sigma}_{n+1}) - \mathbb{E}_x \min(\tau, \overline{\sigma}_n)) \right]$$

$$= \sup_{x \in \gamma} \left[ \mathbb{E}_x \min(\tau, \overline{\sigma}_1) + \sum_{n=1}^\infty \mathbb{E}_x (\chi_{\{\tau > \overline{\sigma}_n\}} \mathbb{E}_{X_{\overline{\sigma}_n}^\varepsilon} \min(\tau, \overline{\sigma}_1)) \right]$$

$$\leq 2\varepsilon^{1/2-\varkappa} \sum_{n=0}^\infty (1 - \tfrac{1}{2}\varepsilon^{1/2-\alpha})^n = 4\varepsilon^{\alpha-\varkappa}.$$

This completes the proof of the lemma.   □

PROOF OF LEMMA 2.5. For a fixed $\widetilde{h}$, let $\widetilde{\tau}$ be the first time when the process $X_t^\varepsilon$ reaches the set $\gamma(\widetilde{h})$. The probability of the event that the process which starts at $x \in \overline{\gamma}$ reaches $\gamma(\widetilde{h})$ before reaching $\gamma$ is estimated in the proof of Lemma 4.4 of [10]. Namely, from formula (35) of [10], it follows that there is a positive $\widetilde{h}$ such that

$$\inf_{x \in \overline{\gamma}} \mathbb{P}_x(\widetilde{\tau} < \sigma) \geq \frac{\varepsilon^\alpha}{2\widetilde{h}}$$

for all $\varepsilon$ which are sufficiently small. Let us now fix a value of $\widetilde{h}$ for which this inequality holds and examine the process $X_t^\varepsilon$ which starts at $x \in \gamma(\widetilde{h})$. By the second part of Lemma 3.1, there is a constant $c(\widetilde{h}) > 0$ such that

$$\sup_{x \in \gamma(\widetilde{h})} \mathbb{E}_x e^{-\sigma} \leq 1 - c(\widetilde{h})$$

for all sufficiently small $\varepsilon$. Due to the Markov property of the process,

$$\sup_{x \in \overline{\gamma}} \mathbb{E}_x e^{-\sigma} \leq 1 - \inf_{x \in \overline{\gamma}} \mathbb{P}_x(\widetilde{\tau} < \sigma) + \inf_{x \in \overline{\gamma}} \mathbb{P}_x(\widetilde{\tau} < \sigma) \sup_{x \in \gamma(\widetilde{h})} \mathbb{E}_x e^{-\sigma}$$

$$\leq 1 - \frac{c(\widetilde{h})\varepsilon^\alpha}{2\widetilde{h}} = 1 - r\varepsilon^\alpha.$$

This completes the proof of the lemma. $\square$

PROOF OF LEMMA 2.6. For $h_1 < h_2$, we denote the set $\{x \in U : h_1 \leq h(x) \leq h_2\}$ by $U(h_1, h_2)$. Let us take numbers $r_1, r_2 \in (0, H(x) - H(A))$, which do not depend on $\varepsilon$, such that $r_1 < r_2$. Take numbers $\alpha' < \alpha''$ such that $\alpha < \alpha' < \alpha'' < \frac{1}{2}$. Let $\gamma' = \gamma(\varepsilon^{\alpha'})$ and $\gamma'' = \gamma(\varepsilon^{\alpha''})$. Let $\overline{\sigma}$ be the first time when the process $X_t^\varepsilon$ reaches $\gamma(r_1)$ or $\gamma'$, whichever happens first. Similarly, let $\overline{\overline{\sigma}}$ be the first time when the process $X_t^\varepsilon$ reaches $\gamma(r_2)$ or $\gamma''$, whichever happens first. For $x \in \overline{\gamma}$, using the Markov property of the process, we can write

$$\mathbb{E}_x \left[ f(h(X_\sigma^\varepsilon)) - f(h(X_0^\varepsilon)) - \int_0^\sigma \mathcal{L} f(h(X_s^\varepsilon)) \, ds \right]$$

$$= \mathbb{E}_x \left[ f(h(X_{\overline{\sigma}}^\varepsilon)) - f(h(X_0^\varepsilon)) - \int_0^{\overline{\sigma}} \mathcal{L} f(h(X_s^\varepsilon)) \, ds \right]$$

(16)
$$+ \mathbb{E}_x \Bigg( \chi_{\{h(X_{\overline{\sigma}}^\varepsilon) = r_1\}}$$

$$\times \mathbb{E}_{X_{\overline{\sigma}}^\varepsilon} \left[ f(h(X_\sigma^\varepsilon)) - f(h(X_0^\varepsilon)) \int_0^\sigma \mathcal{L} f(h(X_s^\varepsilon)) \, ds \right] \Bigg)$$

$$+ \mathbb{E}_x \Bigg( \chi_{\{h(X_{\overline{\sigma}}^\varepsilon) = \varepsilon^{\alpha'}\}}$$

$$\times \, \mathbb{E}_{X_\sigma^\varepsilon} \left[ f(h(X_\sigma^\varepsilon)) - f(h(X_0^\varepsilon)) - \int_0^\sigma \mathcal{L} f(h(X_s^\varepsilon)) \, ds \right].$$

The probability $\sup_{x \in \gamma_1} \mathbb{P}_x(h(X_{\overline{\sigma}}^\varepsilon) = r_1)$ is estimated from above by $c(r_1)\varepsilon^\alpha$ for some constant $c(r_1)$ and all sufficiently small $\varepsilon$ (as follows from formula (35) of [10]). By the first part of Lemma 3.1,

$$\lim_{\varepsilon \to 0} \sup_{x \in \gamma(r_1)} \mathbb{E}_x \left[ f(h(X_\sigma^\varepsilon)) - f(h(X_0^\varepsilon)) - \int_0^\sigma \mathcal{L} f(h(X_s^\varepsilon)) \, ds \right] = 0.$$

Therefore, the second term in the right-hand side of (16) is of order $o(\varepsilon^\alpha)$. We also note that

$$\sup_{x \in \gamma'} \mathbb{E}_x \left[ f(h(X_\sigma^\varepsilon)) - f(h(X_0^\varepsilon)) - \int_0^\sigma \mathcal{L} f(h(X_s^\varepsilon)) \, ds \right]$$

$$\leq c \left( |f(\varepsilon^{\alpha'}) - f(0)| + \sup |Lf| \sup_{x \in \gamma'} \mathbb{E}_x \sigma \right).$$

From Lemma 3.2 it follows that $\sup_{x \in \gamma'} \mathbb{E}_x \sigma = O(\varepsilon^{\alpha'}) = o(\varepsilon^\alpha)$ as $\varepsilon \to 0$. Therefore, the third term on the right-hand side of (16) is also of order $o(\varepsilon^\alpha)$. In order to prove (10), it remains to show that

$$(17) \qquad \sup_{x \in \overline{\gamma}} \left| \mathbb{E}_x \left[ f(h(X_{\overline{\sigma}}^\varepsilon)) - f(h(X_0^\varepsilon)) - \int_0^{\overline{\sigma}} \mathcal{L} f(h(X_s^\varepsilon)) \, ds \right] \right| = o(\varepsilon^\alpha).$$

Let $x_t$ be the deterministic process

$$dx_t = \frac{1}{\varepsilon} v(x_t) \, dt$$

and let $T = T(x)$ be the time it takes the process $x_t$, starting at $x$, to make one rotation along the level set, $T(x) = \inf\{t > 0 : x_t = x\}$.

We shall need several facts about the behavior of the processes $x_t$ and $X_t^\varepsilon$ with an initial point in $U(\varepsilon^{\alpha'}, r_1)$. Let us formulate them here as a separate lemma and then continue with the proof of Lemma 2.6.

LEMMA 3.3. (a) *There are positive constants $c_1$ and $c_2$ such that $c_1 \varepsilon \leq T(x) \leq c_2 \varepsilon |\ln \varepsilon|$ for all sufficiently small $\varepsilon$ and all $x \in U(\varepsilon^{\alpha'}, r_1)$.*

(b) *For any $\delta' > 0$, $R > 0$ and all sufficiently small $\varepsilon$, we have*

$$\mathbb{P}_x \left( \sup_{t \leq T(x)} |h(X_t^\varepsilon) - h(x_t)| > \varepsilon^{1/2 - \delta'} \right) < \varepsilon^R \qquad \text{for all } x \in U(\varepsilon^{\alpha'}, r_1).$$

(c) *For any $\delta' > 0$, $R > 0$ and all sufficiently small $\varepsilon$, we have*

$$\mathbb{P}_x \left( \sup_{t \leq T(x)} |X_t^\varepsilon - x_t| > \varepsilon^{1/2 - \alpha' - \delta'} \right) < \varepsilon^R \qquad \text{for all } x \in U(\varepsilon^{\alpha'}, r_1).$$

The first statement immediately follows from the fact that the critical points of $H$ are nondegenerate. The parts (b) and (c) are proved in the same way as the corresponding statements in Lemma 5.3 below [see formulas (58) and (57)]. Let us proceed with the proof of Lemma 2.6.

We inductively define a sequence of stopping times $T_n$ as follows:

$$T_0 = 0, \qquad T_1 = \min(\overline{\sigma}, T(x)), \qquad T_{n+1} = T_n + T_1(X^\varepsilon_{T_n}).$$

Let replace the time interval $[0, \overline{\sigma}]$ by a larger one, $[0, \widetilde{\sigma}]$, where $\widetilde{\sigma}$ is the first of the stopping times $T_n$ which is greater than or equal to $\overline{\sigma}$, that is,

$$\widetilde{\sigma} = \min_{n\,:\,T_n \geq \overline{\sigma}} T_n.$$

We would like to replace $\overline{\sigma}$ by $\widetilde{\sigma}$ in (17). Note that $\widetilde{\sigma} - \overline{\sigma} \leq c_2 \varepsilon |\ln \varepsilon|$ by part (a) of Lemma 3.3. Parts (a) and (b) of Lemma 3.3 easily imply a statement which is slightly stronger than part (b) of the lemma. Namely, for any $\delta' > 0$, $R > 0$ and all sufficiently small $\varepsilon$, we have

$$\mathbb{P}_x\left(\sup_{t \leq c_2 \varepsilon |\ln \varepsilon|} |h(X^\varepsilon_t) - h(x_t)| > \varepsilon^{1/2 - \delta'}\right) < \varepsilon^R \qquad \text{for all } x \in U(\varepsilon^{\alpha'}, r_1).$$

Therefore,

$$\sup_{x \in \overline{\gamma}} |\mathbb{E}_x[f(h(X^\varepsilon_{\widetilde{\sigma}})) - f(h(X^\varepsilon_{\overline{\sigma}}))]|$$

$$= \sup_{x \in \overline{\gamma}} \left|\mathbb{E}_x \mathbb{E}_{X^\varepsilon_{\overline{\sigma}}} \sup_{0 < t \leq c_2 \varepsilon |\ln \varepsilon|} [f(h(X^\varepsilon_0)) - f(h(X^\varepsilon_t))]\right| = o(\varepsilon^\alpha).$$

Also, by part (a) of Lemma 3.3,

$$\sup_{x \in \overline{\gamma}} \left|\mathbb{E}_x\left[\int_0^{\overline{\sigma}} \mathcal{L}f(h(X^\varepsilon_s))\, ds - \int_0^{\widetilde{\sigma}} \mathcal{L}f(h(X^\varepsilon_s))\, ds\right]\right|$$

$$\leq \sup |Lf| \cdot \sup_{x \in \overline{\gamma}} \mathbb{E}_x(\widetilde{\sigma} - \overline{\sigma}) = o(\varepsilon^\alpha).$$

Therefore, (17) will follow if we show that

$$(18) \qquad \sup_{x \in \overline{\gamma}} \left|\mathbb{E}_x\left[f(h(X^\varepsilon_{\widetilde{\sigma}})) - f(h(X^\varepsilon_0)) - \int_0^{\widetilde{\sigma}} \mathcal{L}f(h(X^\varepsilon_s))\, ds\right]\right| = o(\varepsilon^\alpha).$$

Let

$$(19) \qquad R(\varepsilon) = \sup_{x \in U(\varepsilon^{\alpha'}, r_1)} \frac{|\mathbb{E}_x[f(h(X^\varepsilon_{T_1})) - f(h(X^\varepsilon_0)) - \int_0^{T_1} \mathcal{L}f(h(X^\varepsilon_s))\, ds]|}{\mathbb{E}_x T_1}.$$

We shall prove that $\lim_{\varepsilon \to 0} R(\varepsilon) = 0$. Then, due to the Markov property of the process,

$$\sup_{x \in \overline{\gamma}} \left| \mathbb{E}_x \left[ f(h(X_{\widetilde{\sigma}}^{\varepsilon})) - f(h(X_0^{\varepsilon})) - \int_0^{\widetilde{\sigma}} \mathcal{L} f(h(X_s^{\varepsilon})) \, ds \right] \right|$$

$$\leq \sup_{x \in \overline{\gamma}} \sum_{n=0}^{\infty} \left| \mathbb{E}_x \left( \chi_{\{T_n < \widetilde{\sigma}\}} \right. \right.$$

$$\left. \left. \times \mathbb{E}_{X_{T_n}^{\varepsilon}} \left[ f(h(X_{T_1}^{\varepsilon})) - f(h(X_0^{\varepsilon})) - \int_0^{T_1} \mathcal{L} f(h(X_s^{\varepsilon})) \, ds \right] \right) \right|$$

$$\leq R(\varepsilon) \sup_{x \in \overline{\gamma}} \sum_{n=0}^{\infty} \mathbb{E}_x(\chi_{\{T_n < \widetilde{\sigma}\}} \mathbb{E}_{X_{T_n}^{\varepsilon}} T_1) = R(\varepsilon) \sup_{x \in \overline{\gamma}} \mathbb{E}_x \widetilde{\sigma}.$$

By Lemma 3.2, we have that $\sup_{x \in \overline{\gamma}} \mathbb{E}_x \widetilde{\sigma} = O(\varepsilon^{\alpha})$, and (18) will follow if we show that $\lim_{\varepsilon \to 0} R(\varepsilon) = 0$. From Lemma 3.3, part (b), it follows that

$$(20) \qquad \lim_{\varepsilon \to 0} \sup_{x \in U(\varepsilon^{\alpha'}, r_1)} \frac{T(x)}{\mathbb{E}_x T_1} = 1.$$

Let us study the expression in the numerator of the right-hand side of (19). By Itô's formula,

$$\mathbb{E}_x \left[ f(h(X_{T_1}^{\varepsilon})) - f(h(X_0^{\varepsilon})) - \int_0^{T_1} \mathcal{L} f(h(X_s^{\varepsilon})) \right] ds$$

$$(21) \qquad = \tfrac{1}{2} \mathbb{E}_x \int_0^{T_1} (|\nabla h(X_s^{\varepsilon})|^2 f''(h(X_s^{\varepsilon})) + \Delta h(X_s^{\varepsilon}) f'(h(X_s^{\varepsilon}))) \, ds$$

$$\qquad - \mathbb{E}_x \int_0^{T_1} (a(h(X_s^{\varepsilon})) f''(h(X_s^{\varepsilon})) + b(h(X_s^{\varepsilon})) f'(h(X_s^{\varepsilon}))) \, ds.$$

From the definitions of the coefficients $a(h)$ and $b(h)$, it follows that

$$a(h(x)) = \frac{\int_0^T |\nabla h(x_s)|^2 \, ds}{2T(x)} \quad \text{and} \quad b(h(x)) = \frac{\int_0^T \Delta h(x_s) \, ds}{2T(x)}.$$

We shall need the following estimates on the behavior of the coefficients $a$ and $b$ and the derivatives of the function $f$ near zero.

LEMMA 3.4. *The asymptotic behavior of $a$, $b$ and $f''$ is as follows.*

(a) *There are positive constants $c_1$ and $c_2$ such that*

$$c_1 |\ln h|^{-1} \leq a(h) \leq c_2 |\ln h|^{-1}$$

*for all sufficiently small $h$. Moreover, $a'(h) = o(\frac{1}{h})$ as $h \to 0$.*

(b) *There are positive constants $c_3$ and $c_4$ such that*

$$c_3 \leq b(h) \leq c_4$$

*for all sufficiently small h. Moreover, $b'(h) = o(\frac{1}{h})$ as $h \to 0$.*

(c) *There is a positive constant $c_5$ such that*

$$|f''(h)| \leq c_5 |\ln h|$$

*for all sufficiently small h. Moreover, $|f'''(h)| = o(\frac{|\ln h|^2}{h})$ as $h \to 0$.*

The estimates of $a(h)$ and $b(h)$ and the asymptotics for their derivatives follow from the proof of Lemma 4.5 of [10]. The estimates of $f''(h)$ and $f'''(h)$ are due to the estimates of $a$, $b$ and their derivatives, and to the fact that $a(h)f''(h) + b(h)f'(h) \in \Psi$. Let us proceed with the proof of Lemma 2.6.

We would like to add and then subtract the expression $\mathcal{L}f(h(x))\mathbb{E}_x T_1$ from the right-hand side of (21). First, however, we transform it as follows:

$$\mathcal{L}f(h(x))\mathbb{E}_x T_1$$

$$= \mathbb{E}_x \int_0^{T_1} (a(h(x_s))f''(h(x_s)) + b(h(x_s))f'(h(x_s)))\, ds$$

$$= \mathbb{E}_x \int_0^T (a(h(x_s))f''(h(x_s)) + b(h(x_s))f'(h(x_s)))\, ds + T(x)\psi_1(\varepsilon, x)$$

$$= \tfrac{1}{2}\mathbb{E}_x \int_0^T (|\nabla h(x_s)|^2 f''(h(x_s)) + \Delta h(x_s)f'(x_s))\, ds + T(x)\psi_1(\varepsilon, x)$$

$$= \tfrac{1}{2}\mathbb{E}_x \int_0^{T_1} (|\nabla h(x_s)|^2 f''(h(x_s)) + \Delta h(x_s)f'(x_s))\, ds + T(x)\psi_2(\varepsilon, x).$$

Here, $\psi_1$ and $\psi_2$ are such that

$$\lim_{\varepsilon \to 0} \sup_{x \in U(\varepsilon^{\alpha'}, r_1)} |\psi_i(\varepsilon, x)| = 0, \qquad i = 1, 2.$$

The second and the fourth equalities above are due to Lemma 3.3, part (b) [which implies that $\mathbb{P}_x(T \neq T_1) \leq \varepsilon^R$], Lemma 3.3, part (a) [which bounds $T(x)$ from below] and Lemma 3.4 [which bounds the integrand from above]. Thus, we can rewrite the right-hand side of (21) as

$$\tfrac{1}{2}\mathbb{E}_x \int_0^{T_1} (|\nabla h(X_s^\varepsilon)|^2 f''(h(X_s^\varepsilon)) - |\nabla h(x_s)|^2 f''(h(x_s)))\, ds$$

$$+ \tfrac{1}{2}\mathbb{E}_x \int_0^{T_1} (\Delta h(X_s^\varepsilon)f'(h(X_s^\varepsilon)) - \Delta h(x_s)f'(h(x_s)))\, ds$$

$$+ \mathbb{E}_x \int_0^{T_1} (a(h(x_s))f''(h(x_s)) - a(h(X_s^\varepsilon))f''(h(X_s^\varepsilon)))\, ds$$

$$+ \mathbb{E}_x \int_0^{T_1} (b(h(x_s))f'(h(x_s)) - b(h(X_s^\varepsilon))f'(h(X_s^\varepsilon)))\, ds - T(x)\psi_2(\varepsilon, x).$$

Due to Lemmas 3.3 and 3.4, the absolute value of expectation of each of the four integrals above is estimated by an expression of the form $T(x)\psi(\varepsilon, x)$ such that

$$\lim_{\varepsilon \to 0} \sup_{x \in U(\varepsilon^{\alpha'}, r_1)} |\psi(\varepsilon, x)| = 0.$$

This, together with (20), implies that $\lim_{\varepsilon \to 0} R(\varepsilon) = 0$. This completes the proof of (10).

Let us now prove (11). Let us denote the one-sided derivative of $f(h)$ at $h = 0$ by $f'(0)$. Then,

$$f(h) = f(0) + f'(0)h + o(h) \qquad \text{as } h \to 0$$

and

$$\mathcal{L}f(h) = \frac{1}{k}f'(0) + o(1) \qquad \text{as } h \to 0,$$

where $k$ is the same as in the definition of the operator $\mathcal{L}$. Therefore, we can estimate the left-hand side of (11) as follows:

$$\mathbb{E}_\nu\left[f(h(X_\tau^\varepsilon)) - f(h(X_0^\varepsilon)) - \int_0^\tau \mathcal{L}f(h(X_s^\varepsilon))\, ds\right]$$

$$= f'(0)\varepsilon^\alpha + o(\varepsilon^\alpha) - \frac{1}{k}f'(0)\mathbb{E}_\nu\tau + o(1)\mathbb{E}_\nu\tau = o(\varepsilon^\alpha) \qquad \text{as } \varepsilon \to 0.$$

Here, we used the facts that $0 \le h(X_s^\varepsilon) \le \varepsilon^\alpha$ for $0 \le s \le \tau$, that $\mathbb{E}_\nu\tau = k_2\varepsilon^\alpha(1 + o(1))$ as $\varepsilon \to 0$, where $k_2$ is the same as in Lemma 2.4, and that $k_2 = k$. It remains to prove the last statement of Lemma 2.6.

From Lemma 3.1, part (a), it follows that

$$\sup_{x \in U \setminus V^\varepsilon} \left|\mathbb{E}_x\left[f(h(X_\tau^\varepsilon)) - f(h(X_0^\varepsilon)) - \int_0^\tau \mathcal{L}f(h(X_s^\varepsilon))\, ds\right]\right| \to 0 \qquad \text{as } \varepsilon \to 0.$$

Using Lemma 2.2 and arguments similar to those used in the proof of Lemma 2.4, it is easily seen that

$$\sup_{x \in \mathcal{E} \cup V^\varepsilon} \mathbb{E}_x\tau \to 0 \qquad \text{as } \varepsilon \to 0.$$

Therefore,

$$\sup_{x \in \mathcal{E} \cup V^\varepsilon} \left|\mathbb{E}_x\left[f(h(X_\tau^\varepsilon)) - f(h(X_0^\varepsilon)) - \int_0^\tau \mathcal{L}f(h(X_s^\varepsilon))\, ds\right]\right| \to 0 \qquad \text{as } \varepsilon \to 0.$$

This completes the proof of Lemma 2.6.   $\square$

**4. Estimate on the time of exit from the ergodic component.** Instead of $X_t^\varepsilon$, it will be convenient to consider the same process slowed down by a factor of $\varepsilon$. Thus, we define the following process:

$$\tag{22} d\widetilde{X}_t^\varepsilon = v(\widetilde{X}_t^\varepsilon)\,dt + \sqrt{\varepsilon}\,dW_t, \qquad \widetilde{X}_t^\varepsilon \in \mathbb{T}^2.$$

The first time when the process $\widetilde{X}_t^\varepsilon$ reaches the set $\gamma(\varepsilon^{1/2})$ will be denoted by $\widetilde{\tau}$. In the new notation, Lemma 2.2 can be reformulated as follows.

LEMMA 4.1. *For any $\varkappa > 0$, there exists some $\varepsilon_0 > 0$ such that $\mathbb{E}_x\widetilde{\tau} \leq \varepsilon^{-1/2-\varkappa}$ for $\varepsilon \leq \varepsilon_0$, for all $x \in \mathrm{Cl}(\mathcal{E})$.*

The proof of Lemma 4.1 will rely on a number of technical lemmas. We shall state all of these below and will provide their proofs in the next section. Let us first, however, give an intuitive explanation of the result. If the distance from the original point to the periodic component is strictly positive, then it will take time of order $\varepsilon^{-1}$ for the process $\sqrt{\varepsilon}\,dW_t$ to reach the set $\gamma(\varepsilon^{1/2})$. What helps us to get an extra factor $\varepsilon^{1/2-\varkappa}$ is the presence of the ergodic flow $v$ in the right-hand side of (22). To explain our strategy, take a curve $\Gamma$ transversal to the unperturbed flow $y_t' = v(y_t)$. Then, (22) is a small perturbation of $v$, so near each crossing of $\Gamma$ by the deterministic orbit, there is a crossing by the orbit of (22) (in fact, there are infinitely many random crossings near each deterministic one, but we consider the first of them). Let $B$ be a point on $\Gamma$ which is carried into the saddle point by the unperturbed flow. We want to see how long it takes for the random orbit to get inside the $\sqrt{\varepsilon}$-neighborhood of $B$, since then, the noise will help to get inside the periodic component. Let $Y_n$ denote the consecutive crossings of $\Gamma$ by the deterministic orbit. Then, the average distance between $\{Y_n\}_{n=1}^N$ is $1/N$, so we expect to find one of the points $Y_n$ within $O(1/N)$ from $B$. Next, the distance between the deterministic and random orbits is asymptotically Gaussian with variance $\sqrt{N}\varepsilon$. The above two quantities are of the same order if $1/N = \sqrt{N}\varepsilon$, that is, $N = \varepsilon^{-1/3}$. The probability that a Gaussian random variable with variance $\varepsilon^{1/3}$ hits a segment of length $\sqrt{\varepsilon}$ is of order $\sqrt{\varepsilon}/\varepsilon^{1/3} = \varepsilon^{1/6}$. Hence, we expect that the probability of not reaching the $\varepsilon^{1/2}$-neighborhood of $B$ after $\varepsilon^{-1/3}$ crossings does not exceed $1 - c\varepsilon^{1/6}$, so for $n \gg \varepsilon^{-1/3}$, the probability that the random orbit stays inside the ergodic component for $n$ crossings does not exceed $(1 - c\varepsilon^{1/6})^{n\varepsilon^{1/3}}$. This quantity tends to zero if we consider $n \gg \varepsilon^{-1/2}$, and Lemma 4.1 follows from the Markov property of the process.

In the rest of this section, we provide the precise version of the heuristic arguments given above. For various reasons (one of the reasons being that we want to cover a set of rotation numbers of full measure), the actual statements will be slightly worse than the above heuristics, so the result we get will differ from $\mathbb{E}(\widetilde{\tau}) \sim \varepsilon^{-1/2}$, which could be expected if our heuristic arguments were literally correct, by a factor $\varepsilon^{-\varkappa}$, where $\varkappa$ can be taken arbitrarily small.

Let us now proceed with the rigorous arguments. The main result of [1] implies that there is a smooth closed curve $\Gamma$ on the torus which satisfies the following properties.

(a) $\Gamma$ is homeomorphic to a circle and lies in the interior of $\mathcal{E}$.

(b) $\Gamma$ is transversal to the vector field $v$ at every point of $\Gamma$.

(c) Let $y_t$ be the solution of the equation $y_t' = v(y_t)$ with the initial data $y_0 = x \in \mathcal{E}$. Thus, $y_t$ follows the orbit of the unperturbed flow. We denote by $\sigma$ the first positive time when $y_t$ reaches the set $\Gamma$, that is, $\sigma = \sigma(x) = \inf\{t > 0 : y_t \in \Gamma\}$.

There is a unique point $B \in \Gamma$ such that $\sigma(B) = \infty$. For the solution of the equation $y_t' = v(y_t)$ with the initial data $y_0 = B$, we have $\lim_{t\to\infty} y_t = A$, where $A$ is the saddle point on the torus. In other words, the flow line which starts at $B$ enters the saddle point before reaching $\Gamma$. For the rest of the points $x \in \Gamma$, the time $\sigma(x)$ is finite.

(d) We can define a return map $f : \Gamma \setminus \{B\} \to \Gamma$ as follows. For $x \in \Gamma \setminus \{B\}$, we take the solution $y_t$ of the equation $y_t' = v(y_t)$ with the initial data $y_0 = x$. $f(x)$ is then defined to be equal to $y_{\sigma(x)}$.

We can choose coordinates $\theta : \Gamma \to S^1 = [0, 1)$ such that $\theta(B) = 0$ and the return map is conjugate to the rotation by angle $\rho$, that is, $\theta(f(x)) = \theta(x) + \rho\,(\mathrm{mod}\,1)$. In fact, $d\theta = k\,dH$ on $\Gamma \cap D$, where $k$ is a positive constant ($H$ itself is multivalued on the torus, but $dH = -H'_{x_2}\,dx_1 + H'_{x_1}\,dx_2$ is well defined). For the sake of simplicity of notation and without loss of generality, we shall assume that $k = 1$.

Without loss of generality, we may assume that $\Gamma$ is perpendicular to the vector field $v$ in some neighborhood of $B$ (this will be convenient in the proof of Lemma 4.2).

Let $d(x, y)$ be the metric on $\Gamma$ equal to the distance on the circle between $\theta(x)$ and $\theta(y)$. By considering the image of $\Gamma$ under the action of the flow $-v(x)$ for a sufficiently small time, we can obtain another curve $\Gamma'$ which satisfies the same properties as $\Gamma$ and does not intersect it (see Figure 3). We introduce a sequence of stopping times $\tau_n$, where $\tau_0 = 0$ and $\tau_{n+1}$ is the first time, following $\tau_n$, when the process $\widetilde{X}_t^\varepsilon$ reaches $\Gamma$ after first visiting $\Gamma'$, that is,

$$\tau_{n+1} = \inf_{t \geq \tau_n}\{t : \widetilde{X}_t^\varepsilon \in \Gamma, \widetilde{X}_s^\varepsilon \in \Gamma' \text{ for some } \tau_n \leq s < t\}.$$

Let $\widetilde{\tau}_n = \min(\widetilde{\tau}, \tau_n)$. We can define a Markov chain with the state space $\Gamma$ as follows. Let $X_0 = x$ and $X_n = \widetilde{X}_{\tau_n}^\varepsilon$, where $\widetilde{X}_t^\varepsilon$ is the solution of equation (22) starting at $\widetilde{X}_0^\varepsilon = x$. Then, $X_1, X_2, \ldots$ is a Markov chain on $\Gamma$ (or, if $x \in \Gamma$, we can start the Markov chain with $X_0$. We would like to study the successive returns of this chain to the interval $I_\varepsilon = \{x \in \Gamma : d(x, B) \leq \sqrt{\varepsilon}\}$. We are interested in this because we have the following lemma (see Section A.1).

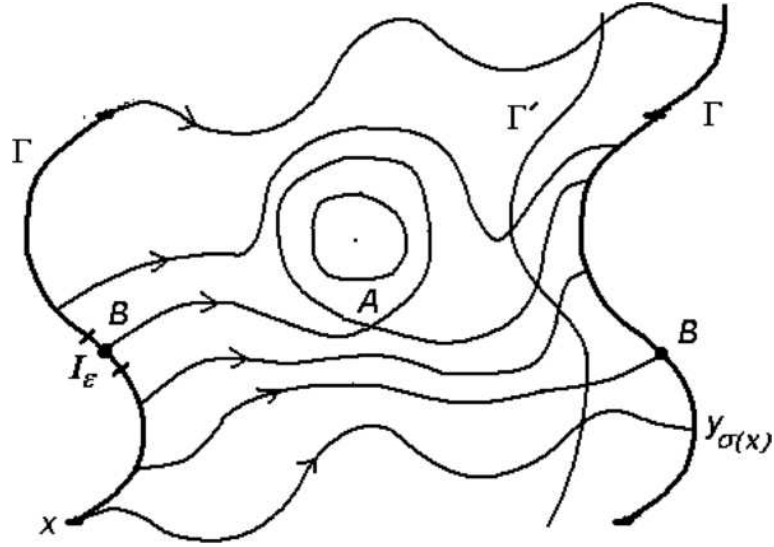

Fig. 3. *Flow on one period.*

LEMMA 4.2. *There exist $c > 0$ and $\varepsilon_0 > 0$ such that $\mathbb{P}_x(\tilde\tau \leq \tau_1) \geq c$ for $\varepsilon \leq \varepsilon_0$, for all $x \in I_\varepsilon$.*

It is easily seen that instead of considering all initial points $x \in \mathcal{E}$ in Lemma 4.1, it is enough to consider $x \in \Gamma$. Indeed, we have the following lemma (see Section A.1).

LEMMA 4.3. *There exist $c > 0$ and $\varepsilon_0 > 0$ such that $\mathbb{E}_x \tilde\tau_1 \leq c|\ln \varepsilon|$ for $\varepsilon \leq \varepsilon_0$, for all $x \in \mathrm{Cl}(\mathcal{E} \cup U_{\varepsilon^{1/2}})$, where $U_{\varepsilon^{1/2}}$ is the part of the periodic component bounded by $\gamma$ and $\gamma(\varepsilon^{1/2})$.*

Let $\tau_I$ be the first of the stopping times $\tau_n$ when $X_n \in I_\varepsilon$, that is, $\tau_I = \min_{n : X_n \in I_\varepsilon} \{\tau_n\}$. We also define $\tilde\tau_I = \min(\tilde\tau, \tau_I)$. We shall prove the following lemma.

LEMMA 4.4. *For any $\varkappa > 0$, there exists some $\varepsilon_0 > 0$ such that $\mathbb{E}_x \tilde\tau_I \leq \varepsilon^{-1/2 - \varkappa}$ for $\varepsilon \leq \varepsilon_0$, for all $x \in \Gamma$.*

Lemmas 4.2–4.4, together with the Markov property of the process $\tilde X_t^\varepsilon$, easily imply Lemma 4.1. Indeed, let $A_n$ be the event that the orbit reaches $\gamma(\varepsilon^{1/2})$ for the first time after exactly $n$ returns to $I_\varepsilon$ and let $\chi_{A_n}$ be the indicator function of this event. Then, $\sup_{x \in I_\varepsilon}(1 - \mathbb{P}_x(A_0)) \leq (1 - c_1)$ and $\sup_{x \in I_\varepsilon} \mathbb{P}_x(A_n) \leq (1 - c_1)^n$, where $c_1$ is the constant from Lemma 4.2. From the relation

$$\sup_{x \in I_\varepsilon} \mathbb{E}_x(\widetilde{\tau} \chi_{A_n}) \leq \sup_{x \in I_\varepsilon} \mathbb{E}_x\left((1 - \chi_{A_0})\left[\widetilde{\tau}_I \sup_{x' \in I_\varepsilon} \mathbb{P}_{x'}(A_{n-1}) + \sup_{x' \in I_\varepsilon} \mathbb{E}_{x'}(\widetilde{\tau} \chi_{A_{n-1}})\right]\right),$$

it follows by induction on $n$ that $\sup_{x \in I_\varepsilon} \mathbb{E}_x(\widetilde{\tau} \chi_{A_n}) \leq (n+1)\varepsilon^{-1/2-\varkappa}(1 - c_1)^{n-1}$. Therefore,

$$\sup_{x \in \mathrm{Cl}(\mathcal{E})} \mathbb{E}_x \widetilde{\tau} \leq \sup_{x \in \mathrm{Cl}(\mathcal{E})} \mathbb{E}_x \widetilde{\tau}_1 + \sup_{x \in \Gamma} \mathbb{E}_x \widetilde{\tau}_I + \sum_{n=0}^{\infty} \sup_{x \in I_\varepsilon} \mathbb{E}_x(\widetilde{\tau} \chi_{A_n})$$

$$\leq c_2 |\ln \varepsilon| + \varepsilon^{-1/2-\varkappa} + \sum_{n=0}^{\infty} (n+1)\varepsilon^{-1/2-\varkappa}(1 - c_1)^{n-1}$$

$$\leq c_3 \varepsilon^{-1/2-\varkappa},$$

where $c_2$ is the constant from Lemma 4.3 and $c_3$ is some other constant. This proves Lemma 4.1.

Let us reduce Lemma 4.4 to a number of simpler statements. Let us consider the deterministic process $Y_n$ on the state space $\Gamma$, where $Y_0 = x$ is the initial point and $Y_{n+1} = f(Y_n)$. Thus, in $\theta$-coordinates, it is simply a rotation by angle $\rho$. Note that $f(Y_n)$ is not defined if $Y_n = B$, but we can define $Y_{n+1}$ in this case using the $\theta$-coordinates. Define

$$N(x) = \min_{n \geq 0}\{n : \sqrt{n\varepsilon} \geq d(Y_n, B)\},$$

where $Y_0 = x$.

We shall need the following lemma about the process $Y_n$ (see Section A.3) (basically, this is a statement about circle rotations by the angle $\rho$).

LEMMA 4.5. *For any $\delta > 0$, there exits some $\varepsilon_0 > 0$ such that $N(x) \leq \varepsilon^{-1/3-\delta}$ for all $\varepsilon \leq \varepsilon_0$, for all $x \in \Gamma$.*

Let us now show that Lemma 4.4 is implied by the following result proved in Section 5.

LEMMA 4.6. *For any $\delta > 0$, there exists some $\varepsilon_0 > 0$ such that*

$$\mathbb{P}_x(X_{N(x)} \in I_\varepsilon) \geq \varepsilon^{1/6+\delta}$$

*for all $\varepsilon \leq \varepsilon_0$, for all $x \in \Gamma$.*

Take an arbitrary $\delta > 0$. From Lemmas 4.5 and 4.6, it follows that

$$\mathbb{P}_x(X_n \notin I_\varepsilon \text{ for all } n \leq \varepsilon^{-1/2-2\delta}) \leq (1 - \varepsilon^{1/6+\delta})^{[\varepsilon^{-1/6-\delta}]}.$$

Therefore, for all sufficiently small $\varepsilon$ and all $x \in \Gamma$, we have

$$\mathbb{P}_x(X_n \notin I_\varepsilon \text{ for all } n \leq \varepsilon^{-1/2-2\delta}) \leq \tfrac{1}{2}. \tag{23}$$

From Lemma 4.3 and the Markov property of the process, it follows that for some $c_1, c_2 > 0$ and all sufficiently small $\varepsilon$, we have

$$\mathbb{P}_x(\widetilde{\tau}_1 > r|\ln \varepsilon|) \leq c_1 e^{-c_2 r}, \qquad r \geq 1,$$

for all $x \in \Gamma$. Therefore,

$$\mathbb{P}_x(\widetilde{\tau}_{[\varepsilon^{-1/2-2\delta}]} > \varepsilon^{-1/2-3\delta}) \leq \tfrac{1}{4} \tag{24}$$

for all sufficiently small $\varepsilon$ for all $x \in \Gamma$. Combining formulas (23) and (24), we obtain

$$\mathbb{P}_x(\widetilde{\tau}_I > \varepsilon^{-1/2-3\delta}) \leq \tfrac{3}{4}$$

for all sufficiently small $\varepsilon$ for all $x \in \Gamma$. With the help of Lemma 4.3, it is easily seen that this, in fact, holds for all $x \in \mathrm{Cl}(\mathcal{E} \cup U_{\varepsilon^{1/2}})$. Using the Markov property of the process, we can now conclude that $\mathbb{E}_x \widetilde{\tau}_I \leq k\varepsilon^{-1/2-3\delta}$ for some constant $k$, for all sufficiently small $\varepsilon$, for all $x \in \mathrm{Cl}(\mathcal{E} \cup U_{\varepsilon^{1/2}})$. Since $\delta$ was an arbitrary positive number, this implies Lemma 4.4.

## 5. Probability of a close encounter.

Here, we prove Lemma 4.6. The statement of this lemma was motivated by the heuristic consideration that $\theta(X_N) - \theta(Y_N)$ is approximately Gaussian with variance of order $\sqrt{N\varepsilon}$. To give a rigorous proof, we need more precise asymptotics. Itô's formula tells us that

$$
\begin{aligned}
\theta(X_N) - \theta(Y_N) &= \sqrt{\varepsilon} \int_0^{\tau_N} \nabla H(\widetilde{X}_s^\varepsilon)\, dW_s + \frac{\varepsilon}{2} \int_0^{\tau_N} \Delta H(\widetilde{X}_s^\varepsilon)\, ds \\
&= \sqrt{\varepsilon} \int_0^{\tau_N} \nabla H(y_s)\, dW_s \\
&\quad + \sqrt{\varepsilon} \int_0^{\tau_N} [\nabla H(\widetilde{X}_s^\varepsilon) - \nabla H(y_s)]\, dW_s \\
&\quad + \frac{\varepsilon}{2} \int_0^{\tau_N} \Delta H(\widetilde{X}_s^\varepsilon)\, ds.
\end{aligned}
\tag{25}
$$

We will show that the main contribution comes from the first term, which is Gaussian. We expect that the last term is $O(\varepsilon N)$ and the second term is a martingale which grows like $\sqrt{V}$, where $V$ is the quadratic variation. We expect that $V$ grows as follows:

$$V \sim \varepsilon \sum_{j=1}^{N} d(X_j, Y_j) \sim \varepsilon \sum_{j=1}^{N} \sqrt{j}\sqrt{\varepsilon} \sim N^{3/2}\varepsilon^{3/2}.$$

This suggests that $\theta(X_N) - \theta(Y_N)$ is approximately Gaussian with correction which is slightly worse than

$$O([(\varepsilon^{-1/3})^{3/2}\varepsilon^{3/2}]^{1/2}) = O(\sqrt{\varepsilon}).$$

This appears bad because we want our process to visit the interval $I_\varepsilon$ of length $O(\sqrt{\varepsilon})$ with positive probability. We overcome this difficulty as follows. We split $N = n_1 + n_2$, where $n_2 \ll \varepsilon^{-1/3}$. During the first $n_1$ iterations, we aim to hit an interval of size $\sqrt{n_2\varepsilon}$, which can be made larger than the correction to Gaussian behavior if $n_2 \gg 1$. During the last $n_2$ steps, we arrange to hit $I_\varepsilon$. Here, it is important that the contribution to the first integral (main Gaussian part) in the right-hand side of (25) from the last $n_2$ steps is larger than the correction accumulated during the first $n_1$ steps. The problem described above does not appear now since the correction accumulated during the last $n_2$ steps is

$$O([n_2^{3/2}\varepsilon^{3/2}]^{1/2}) \ll \sqrt{\varepsilon}.$$

Below, we present the formal, if a bit tedious, implementation of the plan outlined here.

We will need two lemmas concerning the deterministic process $Y_n$. The proofs are given in Section A.3.

**Lemma 5.1.** *For any $\delta > 0$, there exists some $q_0 > 0$ such that*

$$|p - q\rho| \geq q^{-1-\delta}$$

*for all integers $q > q_0$ and $p$.*

**Lemma 5.2.** (a) *There is a constant $c$ such that for any $x \in \Gamma$ and any $N \geq 0$, we have*

$$\sum_{n=0}^{N} \frac{1}{d(Y_n, B)} \leq c\left[N\ln^2 N \left|\ln\left(\min_{0 \leq n \leq N} d(Y_n, B)\right)\right| + \frac{1}{\min_{0 \leq n \leq N} d(Y_n, B)}\right].$$
(26)

(b) *For any $\delta > 0$, there is a constant $c$ such that if*

$$d(Y_n, B) \geq \sqrt{\varepsilon n}$$
(27)

*for all $n \leq N$, then for all sufficiently small $\varepsilon$, we have*

$$\sum_{n=0}^{N} \frac{\sqrt{n}}{d(Y_n, B)} \leq c\varepsilon^{-\delta}[N^{3/2+\delta} + N^\delta \varepsilon^{-1/2}].$$

Let us take some $0 < \alpha < 1/8$. We write $\Gamma = \Gamma_1 \cup \Gamma_2$, where $\Gamma_1$ and $\Gamma_2$ are defined as follows. For $x \in \Gamma_1$, we have $N(x) > 2\varepsilon^{-\alpha}$ and for $x \in \Gamma_2$, we have $N(x) \leq 2\varepsilon^{-\alpha}$. In the first case, we write $N(x) = n_1 + n_2$ with $n_2 = [\varepsilon^{-\alpha}] \leq N(x)/2$. In the second case, we let $n_1 = 0$, $n_2 = N(x)$.

Let $\|\theta\|$ denote the distance from zero of a point $\theta$ on the unit circle and note that there is a naturally defined operation of addition on the circle. We shall need the following lemma (proved in Section A.1). Here, it is convenient to consider the processes on the plane $\mathbb{R}^2$, where the stream function $H$ is single-valued. Let

$$\mathbf{I} = \int_0^{\min(\sigma, \tau_1)} |\nabla H(\tilde{X}_s^\varepsilon) - \nabla H(y_s)|^2 \, ds + \int_{\min(\sigma, \tau_1)}^{\max(\sigma, \tau_1)} |\nabla H(\cdot)|^2 \, ds,$$

where the argument in the second integral is equal to $y_s$ if $\tau_1 < \sigma$ and to $\tilde{X}_s^\varepsilon$ if $\tau_1 > \sigma$.

Let $x^a, x^b \in \Gamma$. Let $y_t^a$ and $y_t^b$ be the orbits of the unperturbed flow, which start at $x^a$ and $x^b$, respectively, and let $\sigma_1^a$ and $\sigma_1^b$ be the corresponding stopping times. Define

$$\mathbf{J}(x^a, x^b) = \int_0^{\min(\sigma_1^a, \sigma_1^b)} |\nabla H(y_s^a) - \nabla H(y_s^b)|^2 \, ds + \int_{\min(\sigma_1^a, \sigma_1^b)}^{\max(\sigma_1^a, \sigma_1^b)} |\nabla H(\cdot)|^2 \, ds,$$

where the argument in the second integral is equal to $y_s^b$ if $\sigma_1^a < \sigma_1^b$ and to $y_s^b$ otherwise.

LEMMA 5.3. *For any $R > 0$ and $\delta', \delta'' > 0$, there exists some $\varepsilon_0 > 0$ such that for $x \in \Gamma$ and all $\varepsilon < \varepsilon_0$, we have*

$$(28) \qquad \mathbb{P}_x\left(\sup_{s < \max(\sigma, \tau_1)} |\tilde{X}_s^\varepsilon - y_s| > \frac{\varepsilon^{1/2 - \delta''}}{d(x, B)}\right) < \varepsilon^R \qquad \text{if } d(x, B) \geq \varepsilon^{1/2 - \delta'},$$

$$(29) \qquad \mathbb{P}_x(\tau_1 > c(\delta')|\ln \varepsilon|) < \varepsilon^R \qquad \text{for some } c(\delta') \text{ if } d(x, B) \geq \varepsilon^{1/2 - \delta'},$$

$$(30) \qquad \mathbb{P}_x(|H(\tilde{X}_{\tau_1}^\varepsilon) - H(x)| > \varepsilon^{1/2 - \delta''}) < \varepsilon^R,$$

$$(31) \qquad \mathbb{E}_x \mathbf{I} < \frac{\varepsilon^{1/2 - \delta''}}{d(x, B)} \qquad \text{if } d(x, B) \geq \varepsilon^{1/2 - \delta'}.$$

*There is a constant $k$ such that for $x^a, x^b \in \Gamma$, we have*

$$(32) \quad \mathbf{J}(x^a, x^b) \leq \frac{k d(x^a, x^b)}{d(x^a, B)} \qquad \text{if } d(x^a, B) > 0 \text{ and } d(x^a, x^b) \leq \frac{1}{2} d(x^a, B).$$

We split the proof of Lemma 4.6 into two parts as follows.

LEMMA 5.4. (a) *There is a positive constant $k_1$ such that for all sufficiently small $\varepsilon$ and all $x \in \Gamma_1$, we have*

$$\mathbb{P}_x(\|\theta(X_{n_1}) + \rho n_2\| \leq \sqrt{n_2 \varepsilon}) \geq \frac{k_1 \sqrt{n_2}}{\sqrt{n_1}},$$

*where $n_1 = n_1(x)$ and $n_2 = n_2(x)$ were defined above.*

(b) *Let $x \in \Gamma$ and $n_2 \leq 2\varepsilon^{-\alpha}$ be such that $\|\theta(x) + \rho n_2\| \leq \sqrt{n_2 \varepsilon}$. There is a positive constant $k_2$, which does not depend on $x$ and $n_2$, such that for all sufficiently small $\varepsilon$, we have*

$$\mathbb{P}_x(X_{n_2} \in I_\varepsilon) \geq \frac{k_2}{\sqrt{n_2}}.$$

Lemma 5.4 implies that for $x \in \Gamma_1$, we have

$$\mathbb{P}_x(X_{N(x)} \in I_\varepsilon) \geq \mathbb{P}_x(A_a)\mathbb{P}_x(A_b|A_a) \geq \frac{k_1 k_2}{\sqrt{n_1}},$$

where $A_a$ and $A_b$ are events from parts (a) and (b) of Lemma 5.4, respectively. Since $n_1 \leq N(x)$, we get the statement of Lemma 4.6 for $x \in \Gamma_1$ from Lemma 4.5. For $x \in \Gamma_2$, we immediately get the statement of Lemma 4.6 from part (b) of Lemma 5.4.

Next, we prove Lemma 5.4. Consider first the case $x \in \Gamma_1$. Let $\sigma_0 = 0$ and let $\sigma_n$ be the time when the trajectory $y_t$ of the unperturbed flow returns to the set $\Gamma$ for the $n$th time (this occurs at the point $Y_n$). Let $s_n = \sigma_{n+1} - \sigma_n$. Let

$$\xi_n = \sqrt{\varepsilon} \int_0^{\sigma_n} \nabla H(y_s)\, dW_s.$$

$\xi_n$ is then Gaussian with zero mean and variance $V_n = \varepsilon \int_0^{\sigma_n} |\nabla H(y_s)|^2\, ds$. Observe that since $\nabla H$ vanishes at $A$, we have

$$(33) \qquad c_1 \leq \int_0^{\sigma_1} |\nabla H(y_s)|^2\, ds \leq c_2$$

and therefore $c_1 \varepsilon n \leq V_n \leq c_2 \varepsilon n$ for some positive constants $c_1$ and $c_2$. We now construct a coupling between $\theta(X_n) - \theta(Y_n)$ and $\xi_n$. Namely, we claim that there exists a random variable $\tilde{\xi}_{n_1}$ with the same distribution as $\xi_{n_1}$ such that for any $\delta, R > 0$, we have

$$(34) \qquad \mathbb{P}_x(\|\theta(X_{n_1}) - \theta(Y_{n_1}) - \tilde{\xi}_{n_1}\| \geq \varepsilon^{1/2-\delta}) = O(\varepsilon^R).$$

This inequality (with $\delta < \alpha$) implies part (a) of Lemma 5.4. Indeed,

$$(35) \qquad \|\theta(Y_{n_1}) + n_2 \rho\| \leq \sqrt{N(x)\varepsilon} \leq \sqrt{2n_1 \varepsilon},$$

so the event of part (a) of Lemma 5.4 happens if $\tilde{\xi}_{n_1}$ belongs to an interval of length $2(\sqrt{n_2 \varepsilon} - \varepsilon^{1/2-\delta})$ whose center is at most $\sqrt{2n_1 \varepsilon}$ away from zero, and $\|\theta(X_{n_1}) - \theta(Y_{n_1}) - \tilde{\xi}_{n_1}\| \leq \varepsilon^{1/2-\delta}$. Formula (34) and the fact that $\tilde{\xi}_{n_1}$ is Gaussian allow us to conclude that the favorable probability is at least $c\sqrt{n_2\varepsilon}/\sqrt{n_1\varepsilon} - O(\varepsilon^R)$ for some positive constant $c$, which implies part (a) of Lemma 5.4. Thus, we need to prove (34).

We now describe the construction of $\tilde{\xi}_n$. Define $\tilde{W}_t^n$ to be $W_t$ for $\tau_n \leq t \leq \tau_{n+1}$ and a Brownian motion (started at $W_{\tau_{n+1}}$) independent of $W_t$ and of $\tilde{W}_t^k$, $k < n$ for $t \geq \tau_{n+1}$. Let

$$(36) \qquad \tilde{\eta}_n = \sqrt{\varepsilon} \int_{\tau_n}^{\tau_n + s_n} \nabla H(y_{\sigma_n + t - \tau_n}) \, d\tilde{W}_t^n, \qquad \tilde{\xi}_n = \sum_{j=0}^{n-1} \tilde{\eta}_j.$$

We want to compare $\tilde{\eta}_n$ with

$$
\begin{aligned}
(37) \qquad H(X_{n+1}) - H(X_n) &= \sqrt{\varepsilon} \int_{\tau_n}^{\tau_{n+1}} \nabla H(\tilde{X}_t^\varepsilon) \, d\tilde{W}_t^n + \frac{\varepsilon}{2} \int_{\tau_n}^{\tau_{n+1}} \Delta H(\tilde{X}_t^\varepsilon) \, dt \\
&= \tilde{\eta}_n' + \tilde{\eta}_n''.
\end{aligned}
$$

Equation (34) will follow if we show that

$$(38) \qquad \mathbb{P}_x \left( \sum_{j=0}^{n_1-1} (\tilde{\eta}_j - (\tilde{\eta}_j' + \tilde{\eta}_j'')) \geq \varepsilon^{1/2-\delta} \right) = O(\varepsilon^R),$$

where $\delta$ and $R$ are the same as in (34). Fix $0 < \alpha_1 < \alpha_2 < \alpha_3 < \delta$ such that $\delta < \frac{3}{2}\alpha_1$. Recall that we assumed that $\alpha < 1/8$, while $\delta < \alpha$ can be taken arbitrarily small. Thus, without loss of generality, we may assume that

$$(39) \qquad 3\alpha + \alpha_1 + 2\delta < \frac{3}{8}.$$

Call $n$ "good" if $d(Y_n, B) \geq \sqrt{n+1}\varepsilon^{1/2-\alpha_1}$ and "bad" otherwise. This distinction between good and bad $n$, is motivated by the fact that for good $n$, the distance $d(X_n, Y_n)$ is much smaller than $d(Y_n, B)$ with large probability, which will be very convenient in the arguments below. Further, formulas (28), (29) and (31) hold if we select $x = Y_n$ as the initial point and $n$ is good.

By Lemma 4.5 and part (b) of Lemma 5.2, the number of bad $n$ such that $n \leq N(x)$ is $O(\varepsilon^{-\alpha_2})$. By (30), the total contribution to the sum in (38) from $\tilde{\eta}_j' + \tilde{\eta}_j''$ with bad $j$ is bounded by $O(\varepsilon^{1/2-\alpha_3})$ [except on a set of probability $O(\varepsilon^R)$]. The contribution to the sum from $\tilde{\eta}_j$ with bad $j$ is a Gaussian random variable with variance bounded from above by $O(\varepsilon^{1-\alpha_2})$, due to (33). Therefore,

$$\max_{0 \leq n \leq n_1} \mathbb{P}_x \left( \left| \sum_{j=0, j-\text{bad}}^{n} (\tilde{\eta}_j - (\tilde{\eta}_j' + \tilde{\eta}_j'')) \right| \geq \frac{\varepsilon^{1/2-\delta}}{2} \right) = O(\varepsilon^R).$$

Since $R$ was arbitrary, this implies that

$$(40) \qquad \mathbb{P}_x\left(\max_{0 \le n \le n_1}\left|\sum_{j=0, j-\text{bad}}^{n}(\tilde{\eta}_j - (\tilde{\eta}_j' + \tilde{\eta}_j''))\right| \ge \frac{\varepsilon^{1/2-\delta}}{2}\right) = O(\varepsilon^R).$$

The contribution to the sum from $\tilde{\eta}_j''$ with good $j$ is bounded from above in absolute value by

$$O(\varepsilon)\sum_{j=0, j-\text{good}}^{n_1-1}(\tau_{j+1} - \tau_j).$$

Therefore, due to Lemma 4.5 and formula (29), we have

$$\max_{0 \le n \le n_1}\mathbb{P}_x\left(\left|\sum_{j=0, j-\text{good}}^{n-1}\tilde{\eta}_j''\right| \ge \varepsilon^{2/3-\alpha_1}\right) = O(\varepsilon^R).$$

Since $R$ was arbitrary, this implies that

$$(41) \qquad \mathbb{P}_x\left(\max_{0 \le n \le n_1}\left|\sum_{j=0, j-\text{good}}^{n-1}\tilde{\eta}_j''\right| \ge \varepsilon^{2/3-\alpha_1}\right) = O(\varepsilon^R).$$

Therefore, formula (38) will follow if we show that

$$(42) \qquad \mathbb{P}_x\left(\sum_{j=0, j-\text{good}}^{n_1-1}(\tilde{\eta}_j - \tilde{\eta}_j') \ge \varepsilon^{1/2-\delta}\right) = O(\varepsilon^R).$$

Let $\hat{y}_t^n$ denote the deterministic orbit starting from $X_n$ at time $\sigma_n$. Let $\hat{s}_n = \sigma(X_n)$. Define

$$\hat{\eta}_n = \sqrt{\varepsilon}\int_{\tau_n}^{\tau_n+\hat{s}_n}\nabla H(\hat{y}_{\sigma_n+t-\tau_n}^n)\,d\tilde{W}_t^n.$$

We want to estimate $\tilde{\eta}_j - \eta_j'$ by comparing both of them to $\hat{\eta}_j$. Let $\ell$ be the first time when $d(X_n, Y_n) > \frac{1}{2}\sqrt{n}\varepsilon^{1/2-\alpha_1}$. Define

$$\tilde{M}_n = \sum_{j=0, j-\text{good}}^{n}(\tilde{\eta}_j - \hat{\eta}_j)\chi_{\{j<\ell\}}, \qquad M_n' = \sum_{j=0, j-\text{good}}^{n}(\eta_j' - \hat{\eta}_j)\chi_{\{j<\ell\}}.$$

The quadratic variation of the martingale $\tilde{M}_n$ is equal to

$$\tilde{V}_n = \varepsilon\sum_{j=0, j-\text{good}}^{n}\chi_{\{j<\ell\}}\left(\int_{\tau_j}^{\tau_j+\min(s_j, \hat{s}_j)}|\nabla H(y_{\sigma_j+t-\tau_j}) - \nabla H(\hat{y}_{\sigma_j+t-\tau_j}^j)|^2\,dt\right.$$

$$\left. + \int_{\tau_j+\min(s_j, \hat{s}_j)}^{\tau_j+\max(s_j, \hat{s}_j)}|\nabla H(\cdot)|^2\,dt\right)$$

(the argument in the second integral is either $y$ or $\hat{y}$, depending on whether $s_j$ or $\hat{s}_j$ is larger). The quadratic variation of the martingale $M'_n$ is

$$V'_n = \varepsilon \sum_{j=0, j-\text{good}}^{n} \chi_{\{j < \ell\}}$$

$$\times \mathbb{E}_x \Bigg[ \bigg( \int_{\tau_j}^{\min(\tau_{j+1}, \tau_j + \hat{s}_j)} |\nabla H(\widetilde{X}_t^\varepsilon) - \nabla H(\hat{y}_{\sigma_j + t - \tau_j}^j)|^2 \, dt$$

$$+ \int_{\min(\tau_{j+1}, \tau_j + \hat{s}_j)}^{\max(\tau_{j+1}, \tau_j + \hat{s}_j)} |\nabla H(\cdot)|^2 \, dt \bigg) \Big| \mathcal{F}_j \Bigg]$$

(the argument in the second integral depends on whether $\tau_{j+1}$ or $\tau_j + \hat{s}_j$ is larger), where $\mathcal{F}_j$ is the $\sigma$-algebra of events determined by $W_t$ by the time $\tau_j$. By (32), we have the estimate

$$\tag{43} \widetilde{V}_{n_1 - 1} \le \varepsilon \sum_{j=0, j-\text{good}}^{n_1 - 1} k \frac{\sqrt{j} \varepsilon^{1/2 - \alpha_1}}{d(Y_j, B)} \le \varepsilon^{1 - \alpha_2},$$

where the second inequality is due to part (b) of Lemma 5.2. We can estimate $V'_n$ by using (31) instead of (32). Namely,

$$\tag{44} V'_{n_1 - 1} \le \varepsilon \sum_{j=0, j-\text{good}}^{n_1 - 1} \frac{2 \varepsilon^{1/2 - \alpha_1}}{d(Y_j, B)} \le \varepsilon^{1 - \alpha_2}.$$

If $R > 0$ is fixed, then applying the $L^p$-maximal inequality, we get

$$\tag{45} \mathbb{P} \bigg( \max_{0 \le n \le n_1} |\tilde{M}_{n-1}| \ge \varepsilon^{1/2 - \alpha_2} \bigg) \le \frac{c(p) \mathbb{E}((\widetilde{V}_{n_1 - 1})^p)}{\varepsilon^{p(1 - 2\alpha_2)}} = O(\varepsilon^{\alpha_2 p}) = O(\varepsilon^R),$$

if $p$ is large enough. A similar argument shows that

$$\tag{46} \mathbb{P} \bigg( \max_{0 \le n \le n_1} |M'_{n-1}| \ge \varepsilon^{1/2 - \alpha_2} \bigg) = O(\varepsilon^R).$$

Therefore, we obtain

$$\mathbb{P}_x \bigg( \sum_{j=0, j-\text{good}}^{n_1 - 1} (\tilde{\eta}_j - \tilde{\eta}'_j) \ge \varepsilon^{1/2 - \delta} \bigg) \le O(\varepsilon^R) + \mathbb{P}_x(\ell < n_1),$$

so it remains to show that the last term is $O(\varepsilon^R)$. To this end, we shall show that the main contribution to the distance between $X_n$ and $Y_n$ comes from the Gaussian term $\widetilde{\xi}_n$.

Define $n_3 = 64 \varepsilon^{-2(\delta - \alpha_1)}$. We have

$$\mathbb{P}_x(\ell < n_1) = \mathbb{P}_x(\ell < n_3) + \mathbb{P}_x(n_3 \le \ell < n_1).$$

The first term can be estimated as follows:

$$\mathbb{P}_x(\ell < n_3) \leq \mathbb{P}_x\left(\max_{0 \leq n \leq n_3} |H(X_{n+1}) - H(X_n)| \geq \frac{\varepsilon^{1/2 - \alpha_1}}{2n_3}\right)$$

$$= \mathbb{P}_x\left(\max_{0 \leq n \leq n_3} |H(X_{n+1}) - H(X_n)| \geq \frac{\varepsilon^{1/2 - (3\alpha_1 - 2\delta)}}{128}\right) = O(\varepsilon^R),$$

by (30). On the other hand by (36), (37), (40) and (41), we have

$$(47) \qquad d(X_n, Y_n)\chi_{\{\ell \geq n\}} \leq |\tilde{\xi}_n - \tilde{M}_{n-1} - M'_{n-1}| + Z_n,$$

where

$$(48) \qquad \mathbb{P}_x\left(\max_{0 \leq n \leq n_1} |Z_n| > \varepsilon^{1/2 - \delta}\right) = O(\varepsilon^R).$$

It follows from (47) that

$$\mathbb{P}_x(n_3 \leq \ell < n_1) \leq \mathbb{P}\left(\max_{n_3 \leq n \leq n_1} \frac{|M'_{n-1}|}{\sqrt{n}} \geq \frac{1}{8}\varepsilon^{1/2 - \alpha_1}\right)$$

$$+ \mathbb{P}\left(\max_{n_3 \leq n \leq n_1} \frac{|\tilde{M}_{n-1}|}{\sqrt{n}} \geq \frac{1}{8}\varepsilon^{1/2 - \alpha_1}\right)$$

$$+ \mathbb{P}_x\left(\max_{n_3 \leq n \leq n_1} \frac{\tilde{\xi}_n}{\sqrt{n}} \geq \frac{1}{8}\varepsilon^{1/2 - \alpha_1}\right)$$

$$+ \mathbb{P}_x\left(\max_{n_3 \leq n \leq n_1} \frac{|Z_n|}{\sqrt{n}} > \frac{1}{8}\varepsilon^{1/2 - \alpha_1}\right).$$

From the definition of $n_3$, it follows that for $n \geq n_3$, we have $\sqrt{n}\varepsilon^{1/2 - \alpha_1}/8 \geq \varepsilon^{1/2 - \delta}$ and $\sqrt{n}\varepsilon^{1/2 - \alpha_1}/8 \geq \varepsilon^{1/2 - \alpha_2}$. Since $\tilde{\xi}_n$ is Gaussian and its variance is bounded by a constant factor of $\varepsilon n$, due to (33), we have

$$\mathbb{P}\left(|\tilde{\xi}_n| \geq \frac{\sqrt{n}\varepsilon^{1/2 - \alpha_1}}{8}\right) \leq k_1 e^{-k_2 \varepsilon^{-2\alpha_1}}$$

for some positive constants $k_1$ and $k_2$, which implies that

$$(49) \qquad \mathbb{P}\left(\max_{n_3 \leq n \leq n_1} |\tilde{\xi}_n| \geq \frac{\sqrt{n}\varepsilon^{1/2 - \alpha_1}}{8}\right) = O(\varepsilon^R).$$

Hence (45), (46), (48) and (49) give

$$\mathbb{P}_x(n_3 \leq \ell \leq n_1) = O(\varepsilon^R),$$

implying (42), which completes the proof of Lemma 5.4(a).

Let us now discuss the modifications needed to prove part (b) of Lemma 5.4. Assuming that $x \in \Gamma$ and $n_2 \leq 2\varepsilon^{-\alpha}$ are such that $\|\theta(x) + \rho n_2\| \leq \sqrt{n_2 \varepsilon}$, we prove that $d(Y_n, B) \geq \varepsilon^{1/8}$ for all $n < n_2$ if $\varepsilon$ is sufficiently small. Indeed, we would otherwise get

$$|(n_2 - n)\rho - p| \leq \sqrt{n_2 \varepsilon} + \varepsilon^{1/8} \leq 2\varepsilon^{1/8}$$

for some integer $p$ since $\alpha < 1/8$. This contradicts Lemma 5.1 since $n_2 - n \leq 2\varepsilon^{-\alpha}$.

We can now repeat the construction used in the proof of part (a) of Lemma 5.4. We claim that there exist $\delta > 0$ and a random variable $\tilde{\xi}_{n_2}$ with the same distribution as $\xi_{n_2}$ such that for any $R > 0$, we have

$$(50) \qquad \mathbb{P}_x(\|\theta(X_{n_2}) - \theta(Y_{n_2}) - \tilde{\xi}_{n_2}\| \geq \varepsilon^{1/2+\delta}) = O(\varepsilon^R).$$

This inequality implies part (b) of Lemma 5.4 in the same way that (34) implies part (a) of the same lemma.

In order to prove (50), we do not need to separate $j < n_2$ into good and bad [instead, we treat all $j$ as we treated the good ones in the proof of part (a)]. The distance from $Y_n$ to $B$ is now controlled using the inequality $d(Y_n, B) \geq \varepsilon^{1/8}$ for all $n < n_2$. The sums in (43) and (44) can now be estimated by

$$O(\varepsilon n_2^{3/2} \varepsilon^{-1/8} \varepsilon^{1/2-\alpha_1}) = O(\varepsilon^{1+[3/8-(3\alpha+\alpha_1)]}) = O(\varepsilon^{1+2\delta}),$$

due to (39). This inequality gives the improvement of $\varepsilon^{1/2+\delta}$ in (50) compared to $\varepsilon^{1/2-\delta}$ in (34). The rest of the technical details remain the same as in the proof of part (a). □

# APPENDIX

**A.1. Proofs of technical lemmas.** In this section, we prove Lemmas 4.2, 4.3 and 5.3.

We shall need the following simple lemma (see, e.g., [10], where it was proven in the case $n = 1$, the proof for $n > 1$ being similar).

LEMMA A.1. Let $X_t^1$ and $X_t^2$ be the following two diffusion processes on $\mathbb{R}^d$ with infinitely smooth coefficients:

$$dX_t^1 = v(X_t^1)\,dt + a(X_t^1)\,dW_t + \sqrt{\varepsilon}v_1(X_t^1)\,dt + \sqrt{\varepsilon}a_1(X_t^1)\,dW_t,$$
$$dX_t^2 = v(X_t^2)\,dt + a(X_t^2)\,dW_t + \sqrt{\varepsilon}v_2(X_t^2)\,dt + \sqrt{\varepsilon}a_2(X_t^2)\,dW_t,$$

with $X_0^1 = X_0^2$. Suppose that for a certain constant $L$, the following bound on the coefficients holds:

$$|\nabla v^i|, |\nabla a^{ij}|, |v_1^i|, |v_2^i|, |a_1^{ij}|, |a_2^{ij}| \leq L, \qquad i, j = 1, \dots, d,$$

where $i$ and $j$ stand for the vector (matrix) entries of the coefficients. Let $\mu$ be the initial distribution for the processes $X_0^1$ and $X_0^2$. Then, for any positive integer $n$, some constant $K = K(L, n)$ and any $t, \eta > 0$, we have

$$\mathbb{P}_\mu\left(\sup_{0 \leq s \leq t}|X_s^1 - X_s^2| \geq \eta\right) \leq \frac{(e^{Kt} - 1)\varepsilon^n}{\eta^{2n}}.$$

PROOF OF LEMMA 4.2. Let us consider the flow on the plane (so as to make $H$ single-valued) and let $A$ be one of the saddle points (there are countably many saddle points on the plane which correspond to the saddle point on the torus). Without loss of generality, we may assume that $H(A) = 0$. Let $B \in I_\varepsilon$ be the point which is carried to $A$ by the deterministic flow (see the picture below). There are two branches of the separatrix that leave the point $A$. One goes inside the ergodic component, while the other forms the boundary between the periodic and the ergodic components. Let us take a point $R$ on the latter, which is distance $d$ away from point $A$ (we measure the distance along the separatrix in the direction of the flow from $A$ to $R$, and $d$ will be specified below). Let us also take a point $Q$ on the separatrix which is distance $d$ away from point $B$ (here, we measure the distance in the direction of the flow from $Q$ to $B$).

Let $\widetilde{\gamma}$ be the curve which consists of two parts of the separatrix—between $Q$ and $A$, and between $A$ and $R$. In a neighborhood of any point $x \in \widetilde{\gamma}$, $x \neq A$, we can consider the smooth change of coordinates $(x_1, x_2) \to (\varphi, \theta)$, where $\varphi = H/\sqrt{\varepsilon}$ and $\theta$ is defined by the following conditions: $|\nabla \theta| = |\nabla H|$ on $\gamma$, $\nabla \theta \perp \nabla H$ and $\theta$ increases in the direction of the deterministic flow. In fact, we can extend this change of coordinates to the region defined by

$$D = \{(\varphi, \theta) : |\varphi| \leq 2, q \leq \theta \leq r\} \setminus \{(\varphi, \theta) : 0 \leq \varphi \leq 2, \theta = a\},$$

where $\theta(Q) = q < \theta(B) = b < \theta(A) = a < \theta(R) = r$. We needed to make the cut in the region $D$ since the change of coordinates degenerates at $A$ (see Figure 4). Let $\tau_D$ be the first time when a trajectory of the process $\widetilde{X}_t^\varepsilon$ reaches the boundary of $D$ [or, rather, the boundary of the preimage of $D$ in $(x_1, x_2)$ coordinates, which will be also denoted by $D$]. If $d$ is small enough, then this set does not intersect $\Gamma'$ and therefore $\tau_D < \tau_1$.

In $(\varphi, \theta)$ coordinates, the generator of the process $\widetilde{X}_t^\varepsilon$ takes the form

$$(51) \quad M^\varepsilon f = \tfrac{1}{2}(f''_{\varphi\varphi}|\nabla H|^2 + \varepsilon f''_{\theta\theta}|\nabla \theta|^2 + \sqrt{\varepsilon} f'_\varphi \Delta H + \varepsilon f'_\theta \Delta \theta) + f'_\theta |\nabla H||\nabla \theta|.$$

The set $I_\varepsilon$ in $(\varphi, \theta)$ coordinates takes the form $\{(\varphi, \theta) : |\varphi| \leq 1, \theta = b\}$. The part of the curve $\gamma(\varepsilon^{1/2})$ that belongs to $D$ takes the form $\{(\varphi, \theta) : \varphi = 1, a < \theta < r\}$. It is sufficient to show that if the process starts at a point $x$ in the former set, then it reaches the latter set before it reaches the boundary of $D$ with probability at least $c > 0$.

Let us take four points on the $\varphi$ axis: $0 > \varphi_1 > \varphi_2 > \varphi_3 > \varphi_4 > -2$. Let $\delta > 0$ be a small constant which will be specified later. Define three segments inside the domain $D$ as follows: $I_1 = \{(\varphi, \theta) : \varphi_2 \leq \varphi \leq \varphi_3, \theta = a - \delta\}$, $I_2 = \{(\varphi, \theta) : \varphi_1 \leq \varphi \leq \varphi_4, \theta = a + \delta\}$ and $I_3 = \{(\varphi, \theta) : \varphi = 1, a < \theta < r\}$. Let $\tau_{I_1}$ be the first time when the process reaches $I_1$; similarly, $\tau_{I_2}$ and $\tau_{I_3}$ are the first times when the process reaches $I_2$ and $I_3$, respectively. By the

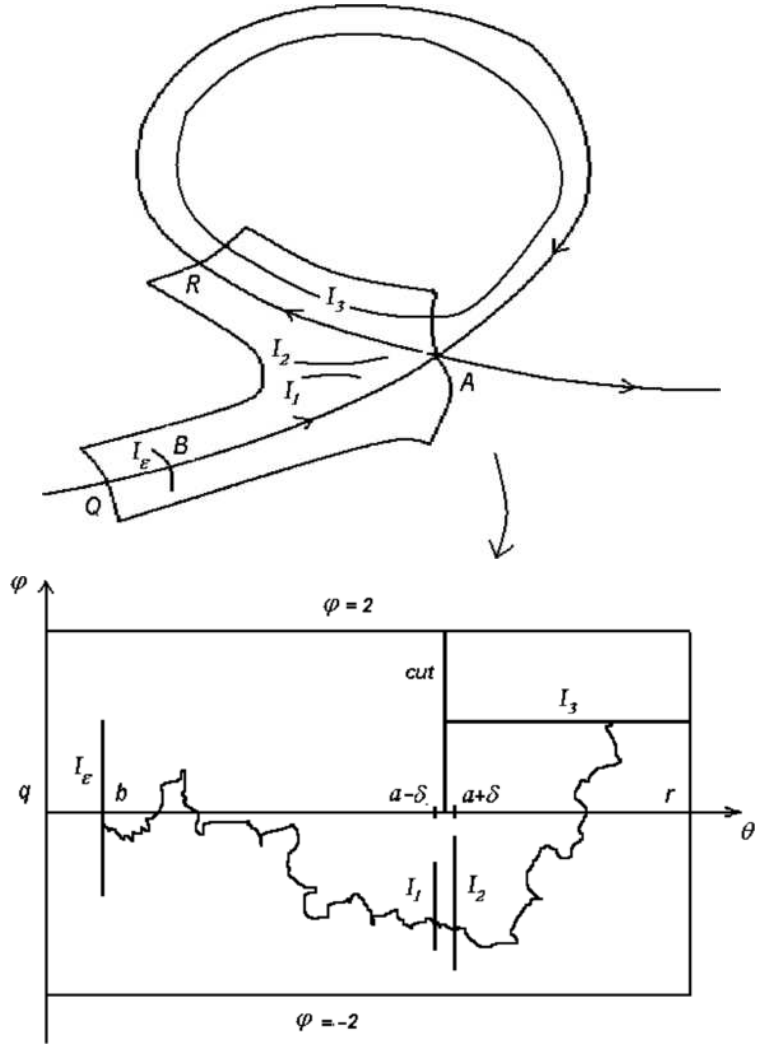

Fig. 4. *Change of variables.*

Markov property of the process, it is sufficient to show that there are positive constants $\delta, c_1, c_2$ and $c_3$ such that

$$
\begin{aligned}
&\inf_{x \in I_\varepsilon} \mathbb{P}_x(\tau_{I_1} < \tau_D) \geq c_1, \\
(52)\quad &\inf_{x \in I_1} \mathbb{P}_x(\tau_{I_2} < \tau_D) \geq c_2 \quad \text{and} \quad \inf_{x \in I_2} \mathbb{P}_x(\tau_{I_3} < \tau_D) \geq c_3.
\end{aligned}
$$

For the first inequality, we define the domain $D_1 = \{(\varphi, \theta) : |\varphi| \leq 2, q \leq \theta \leq a - \frac{\delta}{2}\} \subset D$ (assuming that $\delta$ is fixed). Clearly, $I_\varepsilon, I_1 \subset D_1$. In this domain,

we consider the processes $(\Phi_t, \Theta_t)$ and $(\overline{\Phi}_t, \overline{\Theta}_t)$, which are defined by

$$d\Phi_t = dW_t^\varphi + \frac{\sqrt{\varepsilon}\,\Delta H}{2|\nabla H|^2}\,dt,$$

$$d\Theta_t = \sqrt{\varepsilon}\,\frac{|\nabla\theta|}{|\nabla H|}\,dW_t^\theta + \left(\frac{|\nabla\theta|}{|\nabla H|} + \frac{\varepsilon\Delta\theta}{2|\nabla H|^2}\right)dt$$

and

$$d\overline{\Phi}_t = dW_t^\varphi,$$

$$d\overline{\Theta}_t = dt.$$

The generator of the first process, after multiplying all the coefficients by $|\nabla H|^2$, becomes the operator $M^\varepsilon$ from (51). Therefore, the transition probabilities for this process are the same as for the original process $\widetilde{X}_t^\varepsilon$. We would like to apply Lemma A.1 to the pair of processes $(\Phi_t, \Theta_t)$ and $(\overline{\Phi}_t, \overline{\Theta}_t)$.

Let us follow the process $(\overline{\Phi}_t, \overline{\Theta}_t)$ starting at $x \in I_\varepsilon$ for the time $t_0 = a - \delta/2 - b$. It can easily be seen that the probability that this process reaches $I_1$ before time $t_0$ and before leaving $D_1$ is bounded from below. Moreover, the same is true for all the small perturbations of the trajectories of $(\overline{\Phi}_t, \overline{\Theta}_t)$. More precisely, there are positive constants $\overline{c}_1$ and $\overline{c}_2$ such that for any $x \in I_\varepsilon$, there is an event $\Omega'$, whose probability is at least $\overline{c}_1$, with the property that if $\omega \in \Omega'$, then for any function $(\varphi(t), \theta(t)): [0, t_0] \to D_1$ which satisfies

$$\sup_{t \in [0, t_0]} \|(\varphi(t), \theta(t)) - (\overline{\Phi}_t, \overline{\Theta}_t)(\omega)\| \leq \overline{c}_2,$$

we have $(\varphi(t), \theta(t)) \in D_1$ for all $t \in [0, t_0]$, and $(\varphi(\tau_0), \theta(\tau_0)) \in I_1$ for some $\tau_0 \in [0, t_0]$.

Note that $|\nabla\theta|/|\nabla H| \to 1$ uniformly in $D_1$ as $\varepsilon \to 0$ since we are considering $\delta$ to be fixed for now and the domain $D_1$ is a positive distance away from the saddle point. Therefore, we can consider $(\Phi_t, \Theta_t)$ as a small perturbation of the process $(\overline{\Phi}_t, \overline{\Theta}_t)$ in $D_1$. Let $\Omega''$ be the subset of $\Omega'$ consisting of all $\omega$ for which

$$\sup_{t \in [0, t_0]} \|(\Phi_t, \Theta_t)(\omega) - (\overline{\Phi}_t, \overline{\Theta}_t)(\omega)\| > \overline{c}_2.$$

By Lemma A.1, we can make sure that the probability of $\Omega''$ is less than $\overline{c}_1/2$ for all sufficiently small $\varepsilon$. We thus obtain the first inequality in (52) with $c_1 = \overline{c}_1/2$.

The third inequality in (52) can be proven in exactly the same way if we consider the domain $D_2 = \{(\varphi, \theta) : |\varphi| \leq 2, a + \frac{\delta}{2} \leq \theta \leq r\}$ instead of $D_1$.

Finally, we claim that there is a sufficiently small $\delta$ such that the second inequality holds. Let us consider the intervals $I_1$ and $I_2$ in the original coordinates. By the Morse lemma, there is a smooth change of variables in a neighborhood $O$ of the saddle point $A$ such that in the new variables, the stream function is $H(x_1, x_2) = x_1 x_2$, and the part of $D$ where $\varphi$ was negative now lies in the first quadrant $x_1, x_2 > 0$. In the new variables, the generator of the process $\widetilde{X}_t^\varepsilon$, after a random change of time, becomes $L^\varepsilon f = \varepsilon L_1 f + v_1 \nabla f$, where

$$L_1 f = a_{11}(x_1, x_2) \frac{\partial^2 f}{\partial x_1^2} + a_{12}(x_1, x_2) \frac{\partial^2 f}{\partial x_1 \partial x_2} + a_{22}(x_1, x_2) \frac{\partial^2 f}{\partial x_2^2}$$

$$+ b_1(x_1, x_2) \frac{\partial f}{\partial x_1} + b_2(x_1, x_2) \frac{\partial f}{\partial x_2}$$

is a differential operator with first and second-order terms with bounded coefficients and where $v_1(x_1, x_2) = (-x_1, x_2)$. We shall consider the operator $L^\varepsilon$ in the domain $\mathcal{D}^\varepsilon = O \cap \{x_1 > 0; x_2 > 0; x_1 + x_2 > \varepsilon^{1/3}\}$. We make a further change of variables in $\mathcal{D}^\varepsilon$:

$$(x_1, x_2) \rightarrow (u, v) = \left( \frac{x_1 x_2}{\sqrt{\varepsilon}}, x_2 - x_1 \right).$$

In the new variables, after dividing all the coefficients of the operator by $(x_1 + x_2)$, which amounts to a random change of time for the process, the operator can be written as

$$(53) \qquad\qquad L^\varepsilon f = N^\varepsilon f + \frac{\partial f}{\partial v},$$

where $N^\varepsilon$ is a differential operator with first and second order terms. We claim that all the coefficients of $N^\varepsilon$ are of the form

$$\frac{c_1(x_1, x_2) p_1(x_1, x_2) \sqrt{\varepsilon} + c_2(x_1, x_2) p_2(x_1, x_2) + c_3(x_1, x_2) \varepsilon + c_4(x_1, x_2) \sqrt{\varepsilon}}{x_1 + x_2},$$

where $c_1, \ldots, c_4$ are bounded functions and $p_1$ and $p_2$ are homogeneous first and second degree polynomials, respectively. Indeed, using the expression for $(u, v)$ in terms of $(x_1, x_2)$, we can write

$$\varepsilon a_{12}(x_1, x_2) \frac{\partial^2 f}{\partial x_1 \partial x_2}$$

$$= a_{12}(x_1, x_2) \left( x_1 x_2 \frac{\partial^2 f}{\partial u^2} + (x_2 - x_1) \sqrt{\varepsilon} \frac{\partial^2 f}{\partial u \partial v} + \sqrt{\varepsilon} \frac{\partial f}{\partial u} - \varepsilon \frac{\partial^2 f}{\partial v^2} \right).$$

The other terms of the operator $L_1$ can be treated similarly.

Therefore, all the coefficients of $N^\varepsilon$ can be made arbitrarily small in $\mathcal{D}^\varepsilon$ by first selecting a sufficiently small neighborhood $O$ of the point $A$ and then taking $\varepsilon$ to be sufficiently small.

We can now compare the process whose generator is the operator (53) with the deterministic process with generator $\frac{\partial f}{\partial v}$. If we take $\delta$ to be sufficiently small, then both $I_1$ and $I_2$ are inside $\mathcal{D}^\varepsilon$. Further, the transition time from $I_1$ to $I_2$ for the deterministic process with generator $\frac{\partial f}{\partial v}$ is uniformly bounded in the initial point. The second inequality in (52) can now be deduced from Lemma A.1 in the same way as was the first inequality. $\quad\square$

PROOF OF LEMMA 4.3. The proof of Lemma 4.3 is similar to the arguments which can be found in [7], where the behavior of the process $\widetilde{X}_t^\varepsilon$ in a neighborhood of a saddle point was studied in detail, so we shall only indicate the main steps. As before, we assume that $A$ is one of the saddle points of $H$ on the plane and that $H(A) = 0$. Let $O$ be a small neighborhood of $A$. There are two branches of the separatrix which enter the point $A$ and two which leave it. If $O$ is small enough, then, in a neighborhood of each of the branches of the separatrix intersected with $O$, one can make the change of coordinates $(x_1, x_2) \to (H, \theta)$, where $\theta$ is defined as before, separately in the neighborhood of each of the branches.

Let us take four points, $(P_i, Q_i, R_i, S_i)$, $1 \leq i \leq 4$, on each of the four branches, defined by the conditions $|\theta(P_i) - \theta(A)| = \delta$, $|\theta(Q_i) - \theta(A)| = 2\delta$, $|\theta(R_i) - \theta(A)| = 3\delta$ and $|\theta(S_i) - \theta(A)| = 4\delta$. We can number the branches in such a way that the first and the third branch are the stable directions for the deterministic flow, while the second and the fourth are the unstable directions. Moreover, we can assume that $P_2$ is carried into $P_3$ (and then to $A$) by the flow (see Figure 5).

By the Markov property it suffices to show that there are constants $c, r > 0$ such that for all $x \in \mathrm{Cl}(\mathcal{E} \cup U_{\varepsilon^{1/2}})$ and all sufficiently small $\varepsilon$, we have

$$(54) \qquad \mathbb{P}_x(\widetilde\tau \leq c|\ln\varepsilon|) \geq r.$$

The flow $\widetilde{X}_t^\varepsilon$ can be viewed as a small perturbation of the deterministic flow defined by $y_t' = v(y_t)$, which carries $\Gamma'$ into $\Gamma$ in finite time. This, together with Lemma A.1, implies that

$$(55) \qquad \inf_{x \in \Gamma'} \mathbb{P}_x(\widetilde\tau \leq c_1) \geq r_1$$

for some positive $c_1$ and $r_1$. Let us define four nested neighborhoods of the point $A$ as follows. The neighborhood $\mathcal{P}$ is bounded by eight smooth curves, which are the level sets $|H| = \delta$ and $\theta = \theta(P_i)$, $1 \leq i \leq 4$. Similarly, $\mathcal{Q}, \mathcal{R}$ and $\mathcal{S}$ are bounded by $|H| = 2\delta$ and $\theta = \theta(Q_i)$, $1 \leq i \leq 4$, $|H| = 3\delta$ and $\theta = \theta(R_i)$, $1 \leq i \leq 4$, and $|H| = 4\delta$ and $\theta = \theta(S_i)$, $1 \leq i \leq 4$, respectively.

Let $\tau_{\mathcal{P}}$ be the first time when the process $\widetilde{X}_t^\varepsilon$ reaches $\Gamma'$ or $\mathcal{P}$, whichever happens first. It is clear that the time it takes for the unperturbed process which starts in $\mathrm{Cl}(\mathcal{E} \cup U_{\varepsilon^{1/2}})$ to reach $\Gamma'$ or $\mathcal{P}$ is uniformly bounded in the initial point. This, together with Lemma A.1, easily implies that

$$(56) \qquad \inf_{x \in \mathrm{Cl}(\mathcal{E} \cup U_{\varepsilon^{1/2}})} \mathbb{P}_x(\tau_{\mathcal{P}} \leq c_2) \geq r_2$$

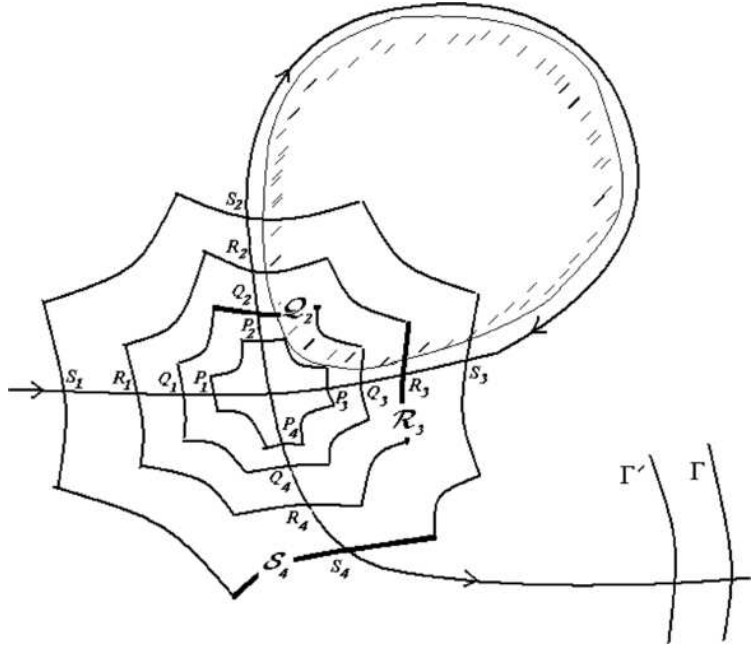

Fig. 5. *Exit from a neighborhood of a saddle point.*

for some positive $c_2$ and $r_2$ and all sufficiently small $\varepsilon$. Due to (55) and (56), we see that it is sufficient to establish (54) for $x \in \mathrm{Cl}(\mathcal{E} \cup U_{\varepsilon^{1/2}}) \cap \mathcal{P}$.

Let $\mathcal{Q}_1$, $\mathcal{Q}_2$, $\mathcal{Q}_3$ and $\mathcal{Q}_4$ be the parts of the boundary of $\mathcal{Q}$ which are given by $\theta = \theta(Q_1)$, $\theta = \theta(Q_2)$, $\theta = \theta(Q_3)$ and $\theta = \theta(Q_4)$, respectively. We shall use similar notation for the parts of the boundary of $\mathcal{R}$ and $\mathcal{S}$. It follows from the arguments in Section 4 of [7] that if the process $\widetilde{X}_t^\varepsilon$ starts at a point $x \in \mathcal{P}$, then it leaves $\mathcal{Q}$ either through $\mathcal{Q}_2$ or through $\mathcal{Q}_4$ with probability which can be made arbitrarily close to one, uniformly in $x \in \mathcal{P}$, by considering sufficiently small $\varepsilon$. Further, the expectation of the time it takes for the process $\widetilde{X}_t^\varepsilon$ to exit $\mathcal{Q}$ is bounded from above by a constant factor of $|\ln \varepsilon|$, as follows from Lemma 4.6 of [7]. Therefore, it is sufficient to establish (54) for $x \in \mathrm{Cl}(\mathcal{E} \cup U_{\varepsilon^{1/2}}) \cap (\mathcal{Q}_2 \cup \mathcal{Q}_4)$.

Any trajectory of the unperturbed flow which starts at $x \in \mathcal{Q}_4$ reaches $\Gamma'$ in finite time and so, with the help of Lemma A.1, we get that

$$\inf_{x \in \mathcal{Q}_4} \mathbb{P}_x(\widetilde{\tau} \leq c_3) \geq r_3$$

for some positive $c_3$ and $r_3$ and all sufficiently small $\varepsilon$. As for initial points in $\mathcal{Q}_2$, the trajectories of the unperturbed flow which start at $x \in \mathcal{Q}_2$ reach $\mathcal{R}_3$ in finite time. Using Lemma A.1 again, we see that it is now sufficient to establish (54) for $x \in \mathrm{Cl}(\mathcal{E} \cup U_{\varepsilon^{1/2}}) \cap \mathcal{R}_3$.

For the process $\widetilde{X}_t^\varepsilon$ which starts at $x \in \mathrm{Cl}(\mathcal{E} \cup U_{\varepsilon^{1/2}}) \cap \mathcal{R}_3$, there is a uniformly positive probability that it exits $\mathcal{S}$ through $\mathcal{S}_4$, as can be seen with the help of the arguments used in the proof of Lemma 4.2. Further, the expectation of the time it takes for the process $\widetilde{X}_t^\varepsilon$ to exit $\mathcal{S}$ is bounded from above by a constant factor of $|\ln \varepsilon|$, as follows from Lemma 4.6 of [7]. Thus, (54) now needs to be established for $x \in \mathcal{S}_4$.

The time it takes for the deterministic process which starts at $x \in \mathcal{S}_4$ to reach $\Gamma'$ is uniformly bounded in $\varepsilon$ and in the initial point. Therefore, we get the desired result by again applying Lemma A.1.   $\square$

PROOF OF LEMMA 5.3.  We first establish two inequalities which are slightly different from (28) and (30), namely

$$(57) \qquad \mathbb{P}_x \left( \sup_{s < \sigma+1} |\widetilde{X}_s^\varepsilon - y_s| > \frac{\varepsilon^{1/2-\delta''}}{d(x,B)} \right) < \varepsilon^R, \qquad \text{if } d(x,B) \geq \varepsilon^{1/2-\delta'},$$

$$(58) \qquad \mathbb{P}_x \left( \sup_{s < \sigma+1} |H(\widetilde{X}_s^\varepsilon) - H(x)| > \varepsilon^{1/2-\delta''} \right) < \varepsilon^R, \qquad \text{if } d(x,B) \geq \varepsilon^{1/2-\delta'}.$$

The proof is based on the use of Lemma A.1. We cannot, however, apply Lemma A.1 to the pair of processes $\widetilde{X}_t^\varepsilon$ and $y_t$ directly, since $\sigma$ grows logarithmically in $\varepsilon$ when $d(x,B) = \varepsilon^{1/2-\delta'}$.

It is clear that $\sigma + 1 \leq c_1 |\ln \varepsilon|$ for some constant $c_1$ which depends on $\delta'$, if $d(x,B) \geq \varepsilon^{1/2-\delta'}$. We claim that the deterministic flow $y_t$ has the following property (this is a minor modification of a similar statement from the proof of Lemma 4.6 of [10], so we do not prove it here).

Let $0 = t_0 < t_1 < t_2 < \cdots < t_m = \sigma + 1$ and $0 < \varkappa < \min(\delta', \delta'')$. Consider a process $\overline{y}_t$ which solves the equation

$$(59) \qquad d\overline{y}_t = v(\overline{y}_t)\, dt$$

on each of the segments $[t_0, t_1), [t_1, t_2), \ldots, [t_{m-1}, t_m]$, with a finite number of jump discontinuities $\lim_{t \to t_i+} \overline{y}_t - \lim_{t \to t_i-} \overline{y}_t = p_i$, $i = 1, \ldots, m-1$. Then, under the conditions

$$y_{t_0} = \overline{y}_{t_0}, \qquad d(y_{t_0}, B) \geq \varepsilon^{1/2-\delta'}, \qquad \sum_{i=1}^{m-1} \|p_i\| < \varepsilon^{1/2-\varkappa},$$

we have

$$(60) \qquad \sup_{0 \leq t \leq t_m} \|y_t - \overline{y}_t\| < c_2 \frac{\sum_{i=1}^{m-1} \|p_i\|}{d(y_{t_0}, B)}$$

for some constant $c_2$ which does not depend on $y_{t_0}$.

Let $K$ be the constant from Lemma A.1 applied to the pair of processes $\widetilde{X}_t^\varepsilon$ and $y_t$, with some $n$ taken to be sufficiently large so that $\varkappa(2n-1) > R$.

Select the points $0 = t_0 < t_1 < t_2 < \cdots < t_m = \sigma + 1$ in such a way that $\frac{\varkappa}{2K}|\ln\varepsilon| \leq |t_{i+1} - t_i| \leq \frac{\varkappa}{K}|\ln\varepsilon|$ for $i = 0, \ldots, m-1$. Since $\sigma + 1 \leq c_1|\ln\varepsilon|$, we have the estimate $m \leq 2c_1K/\varkappa$.

Let $\overline{y}_t^\varepsilon$ be the piecewise continuous process defined by the conditions $\overline{y}_{t_i}^\varepsilon = \widetilde{X}_{t_i}^\varepsilon$ and $d\overline{y}_t^\varepsilon = v(\overline{y}_t^\varepsilon)\,dt$ on $[t_i, t_{i+1}]$, $i = 0, \ldots, m-1$. By Lemma A.1,

$$
\begin{aligned}
(61) \quad \mathbb{P}_x\left(\sum_{i=0}^{m-1} \sup_{t \in [t_i, t_{i+1}]} \|\widetilde{X}_t^\varepsilon - \overline{y}_t^\varepsilon\| > \varepsilon^{1/2-\varkappa}\right) &\leq m\,\frac{e^{\varkappa|\ln\varepsilon|}\varepsilon^n}{(\varepsilon^{1/2-\varkappa}/m)^{2n}} \\
&\leq \left(\frac{2c_1K}{\varkappa}\right)^{1+2n} \varepsilon^{\varkappa(2n-1)}.
\end{aligned}
$$

Due to the continuity of $\widetilde{X}_t^\varepsilon$, formula (61) provides an estimate for the sum of the jumps of the process $\overline{y}_t^\varepsilon$. From (60), it now follows that

$$
(62) \quad \mathbb{P}_x\left(\sup_{0 \leq t \leq \sigma+1} \|y_t - \overline{y}_t^\varepsilon\| > c_2 \frac{\varepsilon^{1/2-\varkappa}}{d(x, B)}\right) \leq \left(\frac{2c_1K}{\varkappa}\right)^{1+2n} \varepsilon^{\varkappa(2n-1)}.
$$

This, together with (61), implies (57) for all sufficiently small $\varepsilon$ since $\varkappa(2n-1) > R$. Since $H(x)$ is constant and $H(\overline{y}_t^\varepsilon)$ is piecewise constant, we have

$$
\begin{aligned}
\mathbb{P}_x&\left(\sup_{0 \leq t \leq \sigma+1} |H(\widetilde{X}_t^\varepsilon) - H(x)| > \varepsilon^{1/2-\delta''}\right) \\
&\leq \mathbb{P}_x\left(\sum_{i=0}^{m-1} \sup_{t \in [t_i, t_{i+1}]} \|\widetilde{X}_t^\varepsilon - \overline{y}_t^\varepsilon\| > \frac{\varepsilon^{1/2-\delta''}}{\sup\|\nabla H\|}\right).
\end{aligned}
$$

This, together with (61) and the condition $\varkappa < \delta''$, implies (58) for all sufficiently small $\varepsilon$.

We take $\delta'' < \delta'$ in (57), and note that if a trajectory of $\widetilde{X}_s^\varepsilon$ stays within a $\frac{\varepsilon^{1/2-\delta''}}{d(x,B)}$-neighborhood of the deterministic trajectory $y_s$ for time $\sigma + 1$, then $\tau_1 < \sigma + 1$. Therefore, (57) and the fact that $\sigma + 1 \leq c_1|\ln\varepsilon|$ together imply (29) with $c = c_1$.

From (57), it follows that $\mathbb{P}_x(\tau_1 < \sigma + 1) \geq 1 - \varepsilon^R$, therefore (57) implies (28).

To prove (30), we consider two cases: $d(x, B) \geq \varepsilon^{1/2-\delta''}$ and $d(x, B) < \varepsilon^{1/2-\delta''}$ [i.e., we set $\delta' = \delta''$ in (58)]. In the first case, the result follows from (58) and the fact that $\mathbb{P}_x(\tau_1 < \sigma + 1) \geq 1 - \varepsilon^R$. To deal with the latter case, we define the set $D$ as the union of all the trajectories of the deterministic flow starting at $x \in \Gamma$ with $d(x, B) \leq \varepsilon^{1/2-\delta''}$ and followed for a positive or negative time until they hit $\Gamma'$. Thus, $\mathrm{Cl}(D)$ contains all the points of $\mathcal{E}$ which are carried too close to the saddle point by the evolution of the deterministic flow starting at $x$ before the time when the trajectory reaches $\Gamma'$.

Observe that our proof of (58) [and consequently of (30)] works not only for initial points $x \in \Gamma$ with $d(x, B) \geq \varepsilon^{1/2 - \delta''}$, but for all points $x \in \mathcal{E} \setminus \mathrm{Cl}(D)$. For $x \in \mathrm{Cl}(D)$, we define $\tau_D$ to be the first time when the process $\widetilde{X}_s^\varepsilon$ reaches the boundary of $U \cup \mathrm{Cl}(D)$. Note that $|H(\widetilde{X}_{\tau_D}^\varepsilon) - H(x)| \leq \varepsilon^{1/2 - \delta''}$ and that (30) is true for $x = \widetilde{X}_{\tau_D}^\varepsilon$. This implies that (30) holds for all $x \in \mathcal{E}$ with $2\varepsilon^{1/2 - \delta''}$, instead of $\varepsilon^{1/2 - \delta''}$, in the left-hand side. Since $\delta''$ was arbitrary, this is equivalent to (30).

Let us now prove (31). The estimate on the expectation of the first term of **I**,

$$\mathbb{E}_x \int_0^{\min(\sigma, \tau_1)} |\nabla H(\tilde{X}_s^\varepsilon) - \nabla H(y_s)|^2 \, ds \leq \overline{c}_1 \frac{\varepsilon^{1/2 - \delta''}}{d(x, B)} \qquad \text{if } d(x, B) \geq \varepsilon^{1/2 - \delta'},$$

is an immediate consequence of (28) (here, $\overline{c}_1$ is some constant). Let $E$ be the event

$$E = \left\{ |\sigma - \tau_1| > k \frac{\varepsilon^{1/2 - \delta''}}{d(x, B)} \right\},$$

where the constant $k$ will be specified below. Then,

$$\mathbb{E}_x \left( (1 - \chi_E) \int_{\min(\sigma, \tau_1)}^{\max(\sigma, \tau_1)} |\nabla H(\cdot)|^2 \, ds \right) \leq \overline{c}_2 \frac{\varepsilon^{1/2 - \delta''}}{d(x, B)} \qquad \text{if } d(x, B) \geq \varepsilon^{1/2 - \delta'}.$$

From (57), it easily follows that for sufficiently large $k$, we have $\mathbb{P}_x(E) \leq \varepsilon^R$ if $d(x, B) \geq \varepsilon^{1/2 - \delta'}$ ($k$ needs to be chosen depending on the minimum of $|\nabla H|$ in a neighborhood of $\Gamma$). Using the fact that the process $\widetilde{X}_{t/\varepsilon}^\varepsilon$ is uniformly (in $\varepsilon$) exponentially mixing, it is not difficult to show that $\mathbb{P}_x(\tau_1 > t/\varepsilon) \leq c e^{-t/c}$ for some positive constant $c$ and all $x \in \mathbb{T}^2$. In particular, $\mathbb{E}_x(\tau_1 \chi_{\{\tau_1 > 1/\varepsilon^2\}}) = o(\varepsilon^R)$ for any positive $R$. Recall that $\sigma \leq c_1 |\ln \varepsilon|$ if $d(x, B) \geq \varepsilon^{1/2 - \delta'}$. Therefore,

$$\mathbb{E}_x \left( \chi_E \int_{\min(\sigma, \tau_1)}^{\max(\sigma, \tau_1)} |\nabla H(\cdot)|^2 \, ds \right) \leq \mathbb{E}_x \left( \chi_E \sigma + \frac{\chi_E}{\varepsilon^2} + \tau_1 \chi_{\{\tau_1 > 1/\varepsilon^2\}} \right)$$

$$\leq \overline{c}_3 \varepsilon^{R-2} \qquad \text{if } d(x, B) \geq \varepsilon^{1/2 - \delta'}.$$

Combining the above estimates, we obtain that (31) holds with an extra constant factor in the right-hand side. Since $\delta''$ was arbitrary, this is equivalent to (31).

Finally, (32) is a statement about deterministic flows. Using the Morse lemma, this statement can be reduced to elementary estimates for a linear system in a neighborhood of the saddle point. $\quad \square$

**A.2. Mixing property for $X_{\tau_n}^\varepsilon$ and $X_{\sigma_n}^\varepsilon$.** This section is devoted to the proof of Lemma 2.3. We need a standard fact concerning the transition density for a diffusion process which starts in the interior of a domain and is stopped when it reaches the boundary. The proof of this lemma is similar to the proof of Theorem 21, I in [12].

LEMMA A.2. *Let $L^\varepsilon(x)$ be a family of differential operators in a connected domain $D$ with a smooth boundary, which are uniformly elliptic in $x \in D$ and $\varepsilon$, and whose coefficients and their first and second derivatives are uniformly bounded in $x \in D$ and $\varepsilon$. Assume that $L^\varepsilon(x)$ are generators for a family of diffusion processes $X_t^\varepsilon$ in the same domain. Let $\sigma$ be the first time when the process $X_t^\varepsilon$ reaches the boundary of $D$. Let $\mu_x^\varepsilon$ be the measure on $\partial D$ defined by $\mu_x^\varepsilon(A) = \mathbb{P}_x(X_\sigma^\varepsilon \in A)$, where $A$ is any Borel subset of $\partial D$. Let $U$ be a domain, whose closure is contained in the interior of $D$.*

*Then, there is a constant $c > 0$, which does not depend on $\varepsilon$, such that $\mu_x^\varepsilon(A) \geq c\lambda(A)$ for all $x \in U$ and all Borel sets $A \subseteq \partial D$, where $\lambda$ is the Lebesgue measure on $\partial D$.*

PROOF OF LEMMA 2.3. Let us prove the result for $\xi_n^2$ (the result for $\xi_n^1$ then follows immediately). As follows from [3], page 197, it is sufficient to prove that there exist a constant $c \geq 0$ and a curve $\widetilde{\gamma} \subseteq \gamma$ (which do not depend on $\varepsilon$) such that $P_2(x, dy) \geq c\lambda(dy)$ for all $x \in \gamma$, $y \in \widetilde{\gamma}$, where $\lambda$ is the Lebesgue measure on $\gamma$.

In the domain between $\gamma$ and $\overline{\gamma}$, we may consider a smooth change of coordinates $(x_1, x_2) \to (\varphi, \theta)$, where $\varphi = H/\sqrt{\varepsilon}$ and $\theta$ is defined by the following conditions: $|\nabla \theta| = |\nabla H|$ on $\gamma$, $\nabla \theta \perp \nabla H$ and $\theta$ increases in the direction of the deterministic flow. We may assume that $\theta \in [0, \int_\gamma |\nabla H| \, dl]$ (with the endpoints of the interval identified) and that the saddle point corresponds to $\theta = 0$. Fix six points $a_1, \ldots, a_6$ on $\gamma$, which satisfy $0 < \theta(a_1) < \cdots < \theta(a_6) < \int_\gamma |\nabla H| \, dl$. We take $\widetilde{\gamma}$ to be the interval between $a_3$ and $a_4$. Let $J$ be the interval [in $(\varphi, \theta)$ coordinates] defined by $1/2 \leq \varphi \leq 1$, $\theta = a_2$.

Any realization of the process, which starts at $x \in \overline{\gamma}$, must pass through the level set $\gamma(\varepsilon^{1/2})$ before hitting $\gamma$. It is easy to show (see, e.g., the arguments in Lemma 3.1 of [10]) that a realization which starts at $x \in \gamma(\varepsilon^{1/2})$ goes through $J$ with probability bounded from below by a positive constant $c_0$ which does not depend on $x$ or $\varepsilon$. Therefore, it is sufficient to show that there is a positive constant $c_1$ such that

$$\mathbb{P}_x(X_\sigma^\varepsilon \in d\theta) \geq c_1 \lambda(d\theta) \tag{63}$$

for all $x \in J$, $\theta \in [a_3, a_4]$ and all sufficiently small $\varepsilon$.

For $\theta_0 \in [a_3, a_4]$ we denote the rectangle $k_1\sqrt{\varepsilon} \leq \varphi \leq 3k_1\sqrt{\varepsilon}$, $\theta_0 - k_2\varepsilon \leq \theta \leq \theta_0 + k_2\varepsilon$ by $\mathcal{R}(\theta_0, \varepsilon, k_1, k_2)$ or simply by $\mathcal{R}$ (see Figure 6). Let $\tau_\mathcal{R}$ be the first time when the process enters this rectangle. We would like to show that there are positive $k_1, k_2$ and $c_2$, which do not depend on $\theta_0$ and $\varepsilon$, such that

$$\mathbb{P}_x(\tau_{\mathcal{R}} < \sigma) \geq c_2\varepsilon \tag{64}$$

if $x \in J$. First, however, let us prove that this inequality implies the statement of the lemma. In $(\varphi, \theta)$ coordinates, the generator of the process $X_t^\varepsilon$ takes the form

$$\begin{aligned}
M^\varepsilon f = & \frac{1}{2\varepsilon}(f_{\varphi\varphi}''|\nabla H|^2 + \varepsilon f_{\theta\theta}''|\nabla\theta|^2 + \sqrt{\varepsilon}f_\varphi'\Delta H + \varepsilon f_\theta'\Delta\theta) \\
& + \frac{1}{\varepsilon}f_\theta'|\nabla H||\nabla\theta|.
\end{aligned} \tag{65}$$

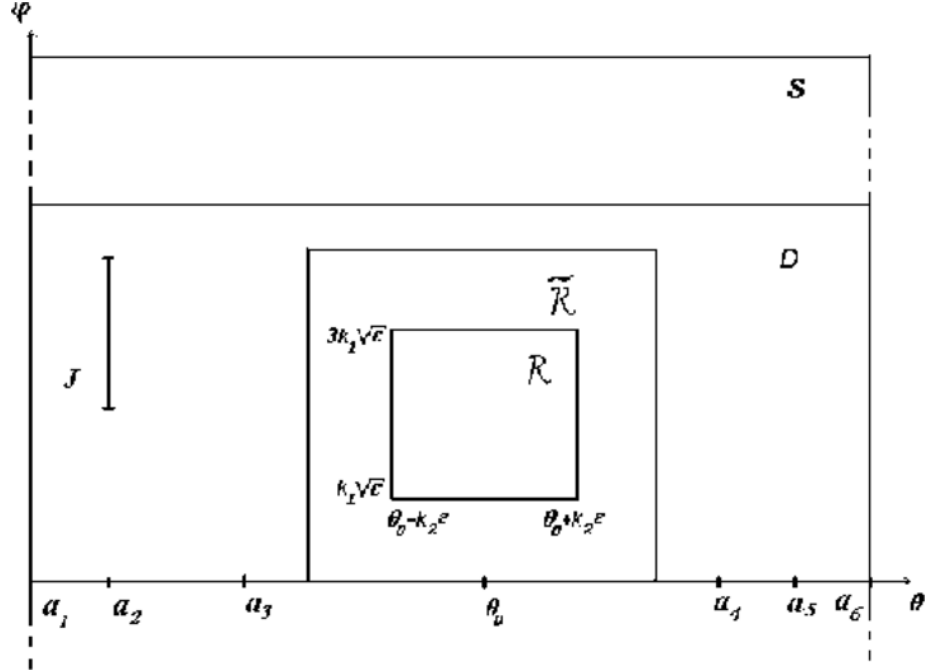

FIG. 6.  *Definition of the rectangles.*

Consider, also, a larger rectangle $0 \leq \varphi \leq 4k_1\sqrt{\varepsilon}$, $\theta_0 - 2k_2\varepsilon \leq \theta \leq \theta_0 + 2k_2\varepsilon$, which will be denoted by $\widetilde{\mathcal{R}}(\theta_0, \varepsilon, k_1, k_2)$ or simply by $\widetilde{\mathcal{R}}$. We can make a further change of variables in $\widetilde{\mathcal{R}}$, namely $(\varphi, \theta) \to (\widetilde{\varphi}, \widetilde{\theta})$, where $\widetilde{\varphi} = \varphi/\sqrt{\varepsilon}$ and $\widetilde{\theta} = (\theta - \theta_0)/\varepsilon$. In the new coordinates, the generator becomes

$$M^\varepsilon f = \frac{1}{2\varepsilon}\left(\frac{1}{\varepsilon}f_{\widetilde{\varphi}\widetilde{\varphi}}''|\nabla H|^2 + \frac{1}{\varepsilon}f_{\widetilde{\theta}\widetilde{\theta}}''|\nabla\theta|^2 + f_{\widetilde{\varphi}}'\Delta H + f_{\widetilde{\theta}}'\Delta\theta\right) + \frac{1}{\varepsilon^2}f_{\widetilde{\theta}}'|\nabla H||\nabla\theta|.$$

After dividing all the coefficients by $|\nabla H|^2/\varepsilon^2$, which amounts to a random change of time for the process, we obtain the new operator for the time-changed process,

$$\overline{M^\varepsilon} f = \frac{1}{2}\left( f''_{\widetilde{\varphi}\widetilde{\varphi}} + f''_{\widehat{\theta}\widehat{\theta}}\frac{|\nabla\theta|^2}{|\nabla H|^2} + \varepsilon f'_{\widetilde{\varphi}}\frac{\Delta H}{|\nabla H|^2} + \varepsilon f'_{\widehat{\theta}}\Delta\theta|\nabla H|^2 \right) + f'_{\widehat{\theta}}\frac{|\nabla\theta|}{|\nabla H|}.$$

We can now apply Lemma A.2 to the process with the generator $\overline{M^\varepsilon}$ in the domain $\widetilde{\mathcal{R}}$ with the initial point $x \in \mathcal{R}$ [note that $\widetilde{\mathcal{R}}$ and $\mathcal{R}$ are fixed sets in $(\widetilde{\varphi},\widetilde{\theta})$ coordinates and that we can smooth out the corners of the rectangle $\widetilde{\mathcal{R}}$ so that it becomes a domain with a smooth boundary]. We conclude that for any $x \in \mathcal{R}$, the measure $\mu_x^\varepsilon(A) = \mathbb{P}_x(X_\sigma^\varepsilon \in A)$, where $A \subseteq \{\widetilde{\varphi} = 0, -1 \le \widetilde{\theta} \le 1\} \subset \partial\widetilde{\mathcal{R}}$, is bounded from below by a measure whose density with respect to $d\widetilde{\theta}$ is equal to a positive constant (which we shall denote by $c_3$).

Since $d\widetilde{\theta} = d\theta/\varepsilon$, the measure $\mu_x^\varepsilon$ (in $\theta$ coordinates) is bounded from below on the interval $[\theta_0 - \varepsilon, \theta_0 + \varepsilon]$ by the measure whose density with respect to $d\theta$ is equal to $c_3/\varepsilon$. Combining this with (64) and using the Markov property of the process, we obtain that (63) holds for all $\theta \in [\theta_0 - \varepsilon, \theta_0 + \varepsilon]$ with $c_1 = c_2 c_3$. Since $\theta_0 \in [a_3, a_4]$ was arbitrary, we see that it remains to prove (64).

We define the following process in $(\varphi, \theta)$ coordinates:

$$d\Phi_t = dW_t^\varphi + \frac{\sqrt{\varepsilon}\,\Delta H}{2|\nabla H|^2}\,dt,$$

$$d\Theta_t = \sqrt{\varepsilon}\frac{|\nabla\theta|}{|\nabla H|}\,dW_t^\theta + \left(\frac{|\nabla\theta|}{|\nabla H|} + \frac{\varepsilon\Delta\theta}{2|\nabla H|^2}\right)dt,$$

where $W_t^\varphi$ and $W_t^\theta$ are independent Brownian motions which can be considered to be defined on different probability spaces $(\Omega^\varphi, \mathcal{F}^\varphi, \mathbb{P}^\varphi)$ and $(\Omega^\theta, \mathcal{F}^\theta, \mathbb{P}^\theta)$. We shall denote the product of these spaces by $(\Omega, \mathcal{F}, \mathbb{P})$. The generator of this process, after multiplying all the coefficients by $|\nabla H|^2/\varepsilon$, becomes the operator $M^\varepsilon$ from (65). Therefore, the transition probabilities for this process are the same as for the original process $X_t^\varepsilon$.

We could also consider the change of variables $(x_1, x_2) \to (\varphi, \theta)$ in an open set which contains the interval $\{H = 0, a_1 \le \theta \le a_6\}$ and does not depend on $\varepsilon$. Thus, the process $(\Phi_t, \Theta_t)$ which starts at $x = (\varphi_0, a_2) \in J$ can be defined until it exits the set $S = \{-\overline{c}/\sqrt{\varepsilon} < \varphi < \overline{c}/\sqrt{\varepsilon}, a_1 < \theta < a_6\}$, where $\overline{c}$ is sufficiently small. Let $\tau_S$ be the first time when the process exits the set $S$. Using Lemma A.1, it is easy to see that $\sup_{x \in J} \mathbb{P}_x(\tau_S \le a_5 - a_2) < \varepsilon^2$ for all sufficiently small $\varepsilon$. Let $\varphi_0$ be fixed and denote the event where $\tau_S > a_5 - a_2$ by $\mathcal{E}_S$.

Let $D$ denote the domain $\{(\varphi, \theta) : 0 < \varphi < 3, a_1 < \theta < a_6\}$ and let the first time when the process leaves this domain be denoted by $\tau_D$. Note that $|\nabla H| = |\nabla \theta|$ for $\varphi = 0$ and that $|\nabla H|$ and $|\nabla \theta|$ are smooth functions of $(H, \theta)$ in $D$, while $H = \sqrt{\varepsilon}\varphi$. Therefore, after rewriting the equations in integral form, we obtain

$$(66) \qquad \Phi_t = \varphi_0 + W_t^\varphi + \sqrt{\varepsilon} \int_0^t g_1^\varepsilon(\Phi_s, \Theta_s)\, ds,$$

$$(67) \qquad \Theta_t = a_2 + \sqrt{\varepsilon}W_t^\theta + \varepsilon \int_0^t g_2^\varepsilon(\Phi_s, \Theta_s)\, dW_s^\theta + t + \sqrt{\varepsilon} \int_0^t g_3^\varepsilon(\Phi_s, \Theta_s)\, ds,$$

where $g_1^\varepsilon$, $g_2^\varepsilon$ and $g_3^\varepsilon$ are bounded in $C^1(D)$ by a constant which does not depend on $\varepsilon$. Note, also, that $g_1^\varepsilon$ is a bounded function in the domain $S$ with a bound which does not depend on $\varepsilon$.

Formulas (66) and (67) show why it is reasonable to expect (64) to hold. They suggest that

$$\Phi_t \approx \varphi_0 + W_t^\varphi,$$

$$\Theta_t \approx a_2 + t + \sqrt{\varepsilon} \int_0^t g_3(\varphi_0 + W_s^\varphi, a_2 + s)\, ds + \sqrt{\varepsilon}W_t^\theta.$$

Let $t_0 = \theta_0 - a_2$. We want $\varphi_0 + W_t^\varphi$ to lie in an interval of size $O(\sqrt{\varepsilon})$, which happens with probability $O(\sqrt{\varepsilon})$. We also want $W_t^\theta$ to cancel $\int_0^t g_3(\varphi_0 + W_s^\varphi, a_2 + s)\, ds$ up to an error of size no larger than $O(\sqrt{\varepsilon})$. The first event, involves only $W_t^\varphi$, while for the second event, we may assume that a realization of $W_t^\varphi$ is fixed and that $W_t^\theta$ is independent of $W_t^\varphi$. So, the required cancellation happens with probability $O(\sqrt{\varepsilon})$. Multiplying those probabilities, we get (64). Let us now proceed with the precise argument.

Consider the time interval $I = [t_0 - \sqrt{\varepsilon}, t_0 + \sqrt{\varepsilon}]$. Let $\mathcal{E}^\varphi \subseteq \Omega^\varphi$ be the event that there is a time $\tau_I$ such that $\tau_I \in I$, $\varphi_0 + W_{\tau_I}^\varphi = 2k_1\sqrt{\varepsilon}$ and $\varphi_0 + W_s^\varphi \in (2k_1\sqrt{\varepsilon}, 2)$ for all $0 \le s < \tau_I$ ($k_1$ will be selected later). Let us note that on the intersections of the events $\mathcal{E}_S \cap (\mathcal{E}^\varphi \times \Omega^\theta)$, we have the following estimates.

(a) Since $g_1^\varepsilon$ is bounded on $S$, there is a constant $\overline{c}_1$ such that

$$\sqrt{\varepsilon} \int_0^{a_5 - a_2} |g_1^\varepsilon(\Phi_s, \Theta_s)|\, ds \le \overline{c}_1 \sqrt{\varepsilon}.$$

Therefore, if we take $k_1 \ge \overline{c}_1$ we obtain $\Phi_{\tau_I} \in [k_1\sqrt{\varepsilon}, 3k_1\sqrt{\varepsilon}]$ and $\Phi_s \in (k_1\sqrt{\varepsilon}, 3)$ for all $0 \le s < \tau_I$. Thus, $(\Phi_s, \Theta_s) \in D$ for $0 \le s < \tau_I$.

(b) Let $\widetilde{\mathcal{E}} \subseteq \Omega$ be the subset of $\mathcal{E}_S \cap (\mathcal{E}^\varphi \times \Omega^\theta)$ where $\sup_{0 \le t \le \tau_I} \varepsilon \int_0^t g_2^\varepsilon(\Phi_s, \Theta_t)\, dW_s^\theta > \sqrt{\varepsilon}$. Since $g_2^\varepsilon$ is bounded on $D$, it is easily seen (e.g., using moment inequalities for martingales) that $\mathbb{P}(\widetilde{\mathcal{E}}) \le \varepsilon^2$ for all sufficiently small $\varepsilon$. Let $\overline{\mathcal{E}}$ be the complement of $\widetilde{\mathcal{E}}$.

(c) Since $g_3^\varepsilon$ is bounded on $D$, there is a constant $\overline{c}_2$ such that

$$\sqrt{\varepsilon} \int_0^{\tau_I} |g_3^\varepsilon(\Phi_s, \Theta_s)|\, ds \le \overline{c}_2 \sqrt{\varepsilon}.$$

Let $\mathcal{E}^\theta \subseteq \Omega^\theta$ be the event that $\sup_{0 \leq s \leq a_5 - a_2} |W_s^\theta| \leq r$, where the constant $r$ will be selected later. On the event $\mathcal{E}_S \cap (\mathcal{E}^\varphi \times \mathcal{E}^\theta) \cap \overline{\mathcal{E}}$, we have

$$\sup_{0 \leq t \leq \tau_I} |\Theta_t - (a_2 + t)| \leq \sqrt{\varepsilon}(r + 1 + \overline{c}_2).$$

Therefore, on this event, we have

$$\Theta_{\tau_I} = a_2 + \sqrt{\varepsilon} W_{\tau_I}^\theta + \varepsilon \int_0^{\tau_I} g_2^\varepsilon(W_s^\varphi, a_2 + s) \, dW_s^\theta$$

$$+ \tau_I + \sqrt{\varepsilon} \int_0^{\tau_I} g_3^\varepsilon(W_s^\varphi, a_2 + s) \, ds + R,$$

where the remainder $R$ satisfies $\mathbb{P}(R > \overline{c}_3 \varepsilon) \leq \varepsilon^2$, the constant $\overline{c}_3$ depending on $r$ [here, we used the fact that $g_2^\varepsilon$ and $g_3^\varepsilon$ are bounded in $C^1(D)$]. For each $\omega^\varphi \in \Omega^\varphi$, the realization $W_s^\varphi$ is fixed and $\tau_I$ is just a constant. Since $g_3^\varepsilon$ is bounded, we have the estimate $\left| \int_0^{\tau_I} g_3^\varepsilon(W_s^\varphi, a_2 + s) \, ds \right| \leq \overline{c}_4$ if $\omega^\varphi \in \mathcal{E}^\varphi$.

As a function of $\omega^\theta \in \Omega^\theta$, the random variable

$$\xi_{\omega^\varphi}(\omega^\theta) = \sqrt{\varepsilon} W_{\tau_I}^\theta + \varepsilon \int_0^{\tau_I} g_2^\varepsilon(W_s^\varphi, a_2 + s) \, dW_s^\theta$$

is Gaussian with the variance equal to $\varepsilon \int_0^{\tau_I} (1 + \sqrt{\varepsilon} g_2^\varepsilon(W_s^\varphi, a_2 + s))^2 \, ds$.

We are interested in the restriction of $\xi_{\omega^\varphi}(\omega^\theta)$ to the set $\mathcal{E}^\theta$. Since $g_2^\varepsilon$ is bounded, it is not difficult to show that for large enough $r$ and all sufficiently small $\varepsilon$, there is a constant $\overline{c}_5$, which does not depend on $\omega^\varphi$ and $\varepsilon$, such that for any Borel set $A \subseteq [-(2 + \overline{c}_4)\sqrt{\varepsilon}, (2 + \overline{c}_4)\sqrt{\varepsilon}]$, we have

$$(68) \qquad \mathbb{P}^\theta(\omega^\theta \in \mathcal{E}^\theta, \xi_{\omega^\varphi}(\omega^\theta) \in A) \geq \frac{\overline{c}_5}{\sqrt{\varepsilon}} \lambda(A),$$

where $\lambda$ is the Lebesgue measure.

Let us collect all the pieces of the proof together. Take $k_1 \geq \overline{c}_1$. With $k_1$ fixed, we have $\mathbb{P}^\varphi(\mathcal{E}^\varphi) \geq q_1 \sqrt{\varepsilon}$ for some constant $q_1$, which does not depend on $\varepsilon$. Let us select $r$ large enough so that (68) holds for all $\omega^\varphi \in \mathcal{E}^\varphi$.

If $(\omega^\varphi, \omega^\theta) \in \mathcal{E}_S \cap (\mathcal{E}^\varphi \cap \mathcal{E}^\theta) \cap \overline{\mathcal{E}}$, then $\Theta_{\tau_I} \in [\theta_0 - k_2 \varepsilon, \theta_0 + k_2 \varepsilon]$, provided that $R < \overline{c}_3 \varepsilon$ and

$$\xi_{\omega^\varphi}(\omega^\theta) \in \left[ \theta_0 - k_2 \varepsilon + \overline{c}_3 \varepsilon - \left( a_2 + \tau_I + \sqrt{\varepsilon} \int_0^{\tau_I} g_3^\varepsilon(W_s^\varphi, a_2 + s) \, ds \right), \right.$$

$$\left. \theta_0 + k_2 \varepsilon - \overline{c}_3 \varepsilon - \left( a_2 + \tau_I + \sqrt{\varepsilon} \int_0^{\tau_I} g_3^\varepsilon(W_s^\varphi, a_2 + s) \, ds \right) \right].$$

Take $k_2 > 2\overline{c}_3$. Then, the interval on the right-hand side of this formula is of length at least $k_2 \varepsilon$ and is centered not further than $(1 + \overline{c}_4)\sqrt{\varepsilon}$ from the origin. Therefore, by (68), for each $\omega^\varphi \in \mathcal{E}^\varphi$ fixed, the $\mathbb{P}^\theta$-probability of this event intersected with $\mathcal{E}^\theta$ is at least $k_2 \overline{c}_5 \sqrt{\varepsilon}$. Since $\mathbb{P}(\mathcal{E}_S \cap \overline{\mathcal{E}} \cap \{R < \overline{c}_3 \varepsilon\}) \geq 1 - 3\varepsilon^2$ and $\mathbb{P}^\varphi(\mathcal{E}^\varphi) \geq q_1 \sqrt{\varepsilon}$, we obtain that

$$\mathbb{P}((\Phi_{\tau_I}, \Theta_{\tau_I}) \in \mathcal{R}, \tau_I < \sigma) \geq q_1 \sqrt{\varepsilon} k_2 \overline{c}_5 \sqrt{\varepsilon} - 3\varepsilon^2 > c_2 \varepsilon,$$

where $c_2 = q_1 k_2 \overline{c}_5/2$ and $\varepsilon$ is sufficiently small. This justifies (64) and thus completes the proof of Lemma 2.3. $\quad\square$

**A.3. Diophantine approximations.** Let us recall some facts about continued fractions (see [9]). Let $\rho = [a_1, \ldots a_n, \ldots]$ be an irrational number between zero and one written as a continued fraction. Let $p_n/q_n = [a_1, a_2 \ldots a_n]$ be the $n$th convergent. Then, $p_{2n}/q_{2n}$ increases to $\rho$ and $p_{2n-1}/q_{2n-1}$ decreases to $\rho$. Also, $q_n = q_{n-1} a_n + q_{n-2}$, so $q_n$ grows at least exponentially. Namely,

$$(69) \qquad q_n \geq f_n \text{ (the } n\text{th Fibonacci number)}.$$

Thus, (2) implies that there is a constant $c$ such that

$$(70) \qquad q_{n+1} \leq c q_n \ln^2 q_n.$$

Let $\|\cdot\|$ denote the distance to the nearest integer. A rational number $p/q$ is called the *best approximation of the second kind* to the number $\rho$ if

$$|p - q\rho| \leq |\tilde{p} - \tilde{q}\rho|$$

for all rational numbers $\tilde{p}/\tilde{q} \neq p/q$ such that $0 < \tilde{q} \leq q$. Thus, $\|q\rho\|$ is minimal among all $\|\tilde{q}\rho\|$ with $0 < \tilde{q} \leq q$. Then (see [9], Section 6), the best approximations of the second kind are exactly the convergents $p_n/q_n$. From here, it follows that

$$(71) \qquad \|q\rho\| < \|q_n\rho\| \qquad \text{implies that } q \geq q_{n+1}.$$

It is known that

$$(72) \qquad \left| \rho - \frac{p_n}{q_n} \right| \leq \frac{1}{q_n q_{n+1}}$$

and

$$(73) \qquad \left| \rho - \frac{p_n}{q_n} \right| \geq \frac{1}{q_n (q_n + q_{n+1})}.$$

Let us show that condition (2) holds for a set of $\rho$ which has Lebesgue measure one. Indeed, let $\mu$ be the measure on $[0,1]$ with the density $((\ln 2)(1 + x))^{-1}$. Then (see [9]),

$$\mu(a_n > n^2) = \mu(a_1 > n^2) = \frac{1}{\ln 2} \ln\left(\frac{n^2 + 2}{n^2 + 1}\right) = \frac{1}{(\ln 2)(n^2 + 1)} + o\left(\frac{1}{n^2 + 1}\right),$$

so, by the Borel–Cantelli lemma, condition (2) holds for $\mu$-almost all $\rho$ and, therefore, also for Lebesgue almost all $\rho$.

We still need to prove Lemmas 4.5, 5.1 and 5.2.

PROOF OF LEMMA 4.5. Let $\hat{N}(x) = \min\{n : d(Y_n, B) \leq \varepsilon^{1/3}\}$. We claim that for any $\delta_2 > 0$, all sufficiently small $\varepsilon$ and all $x \in \Gamma$, we have $\hat{N}(x) \leq$ $\varepsilon^{-(1/3+\delta_2)}$. Applying this claim (with $\delta_2 < \delta$) to $Y_{[\varepsilon^{-1/3}+1]}$ instead of $x$, we see that for all $x$, there exists some $n$ such that $[\varepsilon^{-1/3} + 1] \leq n \leq [\varepsilon^{-1/3} + 1] + \varepsilon^{-(1/3+\delta_2)}$ and $d(Y_n, B) \leq \varepsilon^{1/3}$, thus establishing Lemma 4.5.

To prove the claim, take the smallest $m$ such that $\|q_m\rho\| < \varepsilon^{1/3}$. Suppose that $m$ is even (the case of odd $m$ is almost identical) so that $q_m\rho > p_m$. Consider the interval $J = \{x : 0 \leq \theta(x) \leq \|q_m\rho\|\}$. Let $N_J(x)$ be the first positive time when $Y_n \in J$. Then, $\hat{N} \leq N_J$ so it is enough to show that $N_J(x) \leq \varepsilon^{-(1/3+\delta_2)}$. The maximum of $N_J$ is achieved on $J$ since if $x \notin J$ and $\theta(y) = \theta(x) - \rho$, then $N_J(y) = N_J(x) + 1$.

Note that $N_J$ is piecewise constant on $J$ and its discontinuities are the preimages of the endpoints of $J$ under the first return map to $J$. Our choice of $J$ ensures that there is only one discontinuity inside $J$. Indeed, let $y \in \text{Int}(J)$ be the preimage of the left endpoint. Then, $0 < \theta(y) < \|q_m\rho\|$ and $\theta(y) + s\rho = t$ for some $s, t \in \mathbb{N}$. Then, $\|s\rho\| = \theta(y) < \|q_m\rho\|$. Thus, $s \geq q_{m+1}$ by (71). On the other hand, let $y^*$ be such that

$$\theta(y^*) = \|q_{m+1}\rho\| = p_{m+1} - q_{m+1}\rho \qquad (m+1 \text{ is odd}).$$

If we take $s = q_{m+1}$, then

$$\theta(y^*) + s\rho = p_{m+1} - q_{m+1}\rho + q_{m+1}\rho = p_{m+1},$$

so $y^*$ is the preimage of the left endpoint of $J$.

Likewise, let $y \in \text{Int}(J)$ be the preimage of the right endpoint. Then,

$$\theta(y) + s\rho = q_m\rho - p_m + t, \text{ so}$$

$$(74) \qquad \|s\rho\| = \|q_m\rho\| - \theta(y) < \|q_m\rho\|,$$

$$(75) \qquad \|(s - q_m)\rho\| = \theta(y) < \|q_m\rho\|.$$

Formulas (74) and (71) imply that $s \geq q_{m+1}$ and therefore $s - q_m > 0$. Then, (75) and (71) imply that $s - q_m \geq q_{m+1}$. But,

$$\theta(y^*) + (q_m + q_{m+1})\rho = p_{m+1} + q_m\rho = \|q_m\rho\| + p_{m+1} + p_m.$$

Hence, $y^*$ is also the preimage of the right endpoint. Thus, $N_J$ has only one discontinuity on $J$ and so it takes only two values:

$$N_J(y) = \begin{cases} q_{m+1} + q_m, & \text{if } 0 \leq \theta(y) < \|q_{m+1}\rho\|, \\ q_{m+1}, & \text{if } \|q_{m+1}\rho\| \leq \theta(y) \leq \|q_m\rho\|. \end{cases}$$

From inequality (72), with $m - 1$ instead of $n$, we obtain that $q_m \leq 1/\|q_{m-1}\rho\| \leq \varepsilon^{-1/3}$. Using (70), we obtain $q_m + q_{m+1} \leq \varepsilon^{-(1/3+\delta_2)}$, which completes the proof of the lemma. $\square$

PROOF OF LEMMA 5.1. Formulas (73) and (70) together imply that for some constant $c$, we have

$$(76) \qquad |p_n - q_n\rho| \geq \frac{c}{q_n \ln^2 q_n}.$$

Now, for any $p$ and $q$, we can find $n$ such that $q_n \leq q < q_{n+1}$. We have, therefore,

$$|p - q\rho| \geq |p_{n+1} - q_{n+1}\rho| \geq \frac{c}{q_{n+1} \ln^2 q_{n+1}}$$

$$\geq \frac{c'}{q_n \ln^2 q_n \ln^2(q_n \ln^2 q_n)} \geq \frac{c'}{q \ln^2 q \ln^2(q \ln^2 q)}$$

for some $c'$. This implies the statement of the lemma. $\square$

PROOF OF LEMMA 5.2. Take the smallest $n$ such that $q_n > N$. Recall that the Denjoy–Koksma inequality (see [2], Lemma 1 of Section 3.4) says that if $q_n$ is the denominator of the $n$th convergent, then for any function $\phi$ and any points $x_1, x_2$ on the unit circle,

$$\left| \sum_{j=0}^{q_n-1} [\phi(x_1 + j\rho) - \phi(x_2 + j\rho)] \right| \leq \mathrm{Var}(\phi),$$

where $\mathrm{Var}(\phi)$ denotes the variation of $\phi$. Let $\phi(z) = 1/\|z\|$ if $\|z\| \geq \min_{n \leq N} d(Y_n, B)$ and $\phi(z) = 1/\min_{n \leq N} d(Y_n, B)$ otherwise. Applying the Denjoy–Koksma inequality with $x_1 = \theta(x)$ and $x_2$ such that

$$\sum_{j=0}^{q_n-1} \phi(x_2) = q_n \int_0^1 \phi(z)\,dz,$$

we get

$$\sum_{n=0}^{N} \frac{1}{d(Y_n, B)} \leq \sum_{n=0}^{q_n-1} \frac{1}{d(Y_n, B)} \leq q_n \int_0^1 \phi(z)\,dz + \mathrm{Var}(\phi).$$

Now,

$$\int_0^1 \phi(z)\,dz \sim 2\ln\left(\frac{1}{\min_{n \leq N} d(Y_n, B)}\right), \qquad \mathrm{Var}(\phi) \leq \frac{2}{\min_{n \leq N} d(Y_n, B)}.$$

By (70), we have $q_n \leq cq_{n-1} \ln^2 q_{n-1} \leq cN \ln^2 N$, which proves (a).

To prove (b), let $J_m = [[N^{\delta m/2}], [N^{\delta(m+1)/2}]]$, where $0 \leq m \leq [2/\delta]$. The union of these segments covers the interval $[1, N]$. Since $n \leq [N^{\delta(m+1)/2}]$ on $J_m$, it follows by part (a), applied to the initial point $Y_{[N^{\delta m/2}]}$ instead of $x$, that there is a constant $c$ such that

$$\sum_{n \in J_m} \frac{\sqrt{n}}{d(Y_n, B)}$$

$$\leq c\left[ N^{3\delta(m+1)/4} \ln^2(N^{\delta(m+1)/2}) \left| \ln\left(\min_{n \in J_m} d(Y_n, B)\right) \right| + \frac{N^{\delta(m+1)4}}{\min_{n \in J_m} d(Y_n, B)} \right].$$

Combining this with (27), we get, for some other constant $c$,

$$\sum_{n \in J_m} \frac{\sqrt{n}}{d(Y_n, B)} \le c \left[ N^{3/2+\delta} |\ln(\varepsilon)| + \frac{N^{\delta/4}}{\sqrt{\varepsilon}} \right].$$

Summation over $m$ proves part (b) of the lemma. $\square$

**A.4. The case of several saddle points.** In this section, we discuss the case of several periodic components. First, we assume that each of the domains $U_k$, $k = 1, \ldots, n$, contains a single critical point $M_k$ (a maximum or a minimum of $H$). Let $A_k$, $k = 1, \ldots, n$, be the saddle points of $H$, such that $A_k$ is on the boundary of $U_k$. We denote the boundary of $U_k$ by $\gamma_k$. Let $p_k = \pm \int_{\gamma_k} |\nabla H| \, dl$, where the sign $+$ is taken if $A_k$ is a local minimum for $H$ restricted to $U_k$, and $-$ is taken otherwise.

The phase space of the limiting process is now a graph $\mathbb{G}$ which consists of $n$ edges $I_k$, $k = 1, \ldots, n$ (segments labeled by $k$), where each segment is either $[H(M_k) - H(A_k), 0]$ (if $M_k$ is a minimum) or $[0, H(M_k) - H(A_k)]$ (if $M_k$ is a maximum). All the edges share a common vertex (the origin). Thus, a point on the graph can be determined by specifying $k$ (the number of the edge) and the coordinate on the edge. We define the mapping $h: \mathbb{T}^2 \to \mathbb{G}$ as follows:

$$h(x) = \begin{cases} 0, & \text{if } x \in \mathcal{E} \\ (k, H(x) - H(A)), & \text{if } x \in U_k. \end{cases}$$

We shall use the notation $h_k$ for the coordinate on $I_k$. As in the case of one periodic component, we define the limiting process via its generator $\mathcal{L}$. First, for each $k$, we define the differential operator $L_k f = a_k(h_k) f'' + b_k(h_k) f'$ on the interior of $I_k$, where the coefficients $a_k$ and $b_k$ are given by formulas (3) and (4) [where $\gamma(h_k)$ is defined for each of the periodic components and has the same meaning as in the case of one periodic component]. The domain of $\mathcal{L}$ consists of those functions $f \in C(\mathbb{G})$ which

(a) are twice continuously differentiable in the interior of each of the edges;

(b) have the limits $\lim_{h_k \to 0} L_k f(h_k)$ and $\lim_{h_k \to (H(M_k) - H(A_k))} L_k f(h_k)$ at the endpoints of each of the edges, the value of the limit $q = \lim_{h_k \to 0} L_k f(h_k)$ being the same for all edges;

(c) have the limits $\lim_{h_k \to 0} f'(h_k)$ and

$$\sum_{k=1}^{n} p_k \lim_{h_k \to 0} f'(h_k) = 2 \operatorname{Area}(\mathcal{E}) q.$$

For functions $f$ which satisfy the above three properties, we define $\mathcal{L}f = L_k f$ in the interior of each edge, and as the limit of $L_k f$ at the endpoints of $I_k$.

As in the case of one periodic component, we have the following theorem.

THEOREM 2.    *The measure on $C([0, \infty), \mathbb{G})$ induced by the process $Y_t^\varepsilon = h(X_t^\varepsilon)$ converges weakly to the measure induced by the process with the generator $\mathcal{L}$ with the initial distribution $h(X_0^\varepsilon)$.*

The proof of this theorem requires some modifications to the proof of Theorem 1. We sketch these modifications without providing all of the technical details.

(I) Recall the definition of the Markov chains $\xi_n^1$ and $\xi_n^2$ from Section 2. In the case of several periodic components, the state spaces for these Markov chains will be slightly different. Namely, we replace the curves $\gamma$ and $\overline{\gamma}$ defined in Section 2 by the curves $\gamma = \bigcup_k \gamma_k$ and $\overline{\gamma} = \bigcup_k \overline{\gamma}_k$, where $\overline{\gamma}_k = \{|H_k| = \varepsilon^\alpha\}$. In the proof of Lemma 2.4, $U$ will now stand for the union $U = \bigcup_k U_k$ of the loops. Let $\mu$ and $\nu$ be the invariant measures on $\overline{\gamma}$ and $\gamma$, respectively. Let us study the asymptotics $\mu(\overline{\gamma}_k)$ for different $k$.

Let $\mu_k$ be the normalized restriction of the measure $\mu$ to $\overline{\gamma}_k$, that is,

$$\mu_k(A) = \mu(A)/\mu(\overline{\gamma}_k)$$

for each measurable subset $A$ of $\overline{\gamma}_k$. Instead of (7), we now have

$$(77) \qquad \mathbb{E}_{\mu_k} \sigma = 2\left(\int_{\gamma_k} |\nabla H| \, dl\right)^{-1} \mathrm{Area}(U_k) \varepsilon^\alpha (1 + o(1)) \qquad \text{as } \varepsilon \to 0.$$

Let us prove that

$$(78) \qquad \mu(\overline{\gamma}_k) = \frac{\int_{\gamma_k} |\nabla H| \, dl}{\int_\gamma |\nabla H| \, dl} (1 + o(1)) \qquad \text{as } \varepsilon \to 0.$$

Let $\tau^{(k)}$ be the time when the process $X_t^\varepsilon$ visits $\overline{\gamma}_k$ for the first time. Similarly to (15), we obtain

$$(79) \qquad\qquad \frac{\mathbb{E}_\nu \tau}{\mathbb{E}_\mu \sigma} \sim \frac{\mathrm{Area}(\mathcal{E})}{\mathrm{Area}(U)}$$

and

$$(80) \qquad\qquad \frac{\mathbb{E}_\nu \tau^{(k)}}{\mathbb{E}_{\mu_k} \sigma} \sim \frac{\mathrm{Area}(\mathbb{T}^2 - U_k)}{\mathrm{Area}(U_k)} \qquad \text{as } \varepsilon \to 0.$$

Let $N_k(T)$ be the number of times the process $X_t^\varepsilon$ travels from $\gamma$ to $\overline{\gamma}_k$ before time $T$ and $M(T)$ the number of times the process travels from $\gamma$ to $\overline{\gamma}$ before time $T$. Note that

$$\lim_{T \to \infty} \frac{N_k(T)}{M(T)} = \mu(\overline{\gamma}_k)(1 + o(1)) \qquad \text{as } \varepsilon \to 0.$$

By the Birkhoff ergodic theorem,

$$\lim_{T\to\infty}\frac{T}{N_k(T)}=(\mathbb{E}_\nu\tau^{(k)}+\mathbb{E}_{\mu_k}\sigma)(1+o(1))\qquad\text{as }\varepsilon\to 0$$

and

$$\lim_{T\to\infty}\frac{T}{M(T)}=(\mathbb{E}_\nu\tau+\mathbb{E}_\mu\sigma)(1+o(1))\qquad\text{as }\varepsilon\to 0.$$

Therefore,

$$\lim_{\varepsilon\to 0}\frac{\mu(\overline{\gamma}_k)(\mathbb{E}_\nu\tau^{(k)}+\mathbb{E}_{\mu_k}\sigma)}{\mathbb{E}_\nu\tau+\mathbb{E}_\mu\sigma}=1.$$

Now, (77) and (80) imply that the ratio $\mu(\overline{\gamma}_k)/\int_{\gamma_k}|\nabla H|\,dl$ is asymptotically independent of $k$, thus proving (78).

(II) Let us use (78) to justify (11). Near the origin, we have

$$f(h_k)=f(0)+\lim_{h_k\to 0}f'(h_k)h_k+o(h_k).$$

Let $r_k=1$ if $M_k$ is a maximum and $r_k=-1$ if $M_k$ is a minimum. Observe that since $\mu$ and $\nu$ are invariant, $\nu(X_\tau^\varepsilon\in\overline{\gamma}_k)=\mu(\overline{\gamma}_k)$. Therefore,

$$\mathbb{E}_\nu\left(f(h(X_\tau^\varepsilon))-f(h(X_0^\varepsilon))-\int_0^\tau(Lf)(X_s^\varepsilon)\,ds\right)$$

$$\sim\sum_k r_k\mu(\overline{\gamma}_k)\lim_{h_k\to 0}f'(h_k)-q\mathbb{E}_\nu\tau+o(\varepsilon^\alpha)$$

$$=\varepsilon^\alpha\left[\frac{\sum_k r_k(\int_{\gamma_k}|\nabla H|\,dl)\lim_{h_k\to 0}f'(h_k)}{\int_\gamma|\nabla H|\,dl}-\frac{2q\mathrm{Area}(\mathcal{E})}{\int_\gamma|\nabla H|\,dl}+o(1)\right]=o(\varepsilon^\alpha),$$

where the expression for $\mathbb{E}_\nu\tau$ was obtained from (77), (78) and (79).

(III) Parts (I) and (II) explain the new gluing conditions. The proof of the mixing of $\xi_n^1$ and $\xi_n^2$ also needs to be modified. Namely, we introduce the stopping times $\tau_{\varepsilon,k}=\min\{n\!:\!d(X_n,B_k)\leq\varepsilon^{1/2}\}$, where $B_k$ is the preimage of $A_k$ on $\Gamma$. Also, let $N_k(x)=\min\{n\!:\!d(Y_n,B_k)\leq\sqrt{n\varepsilon}\}$. Lemma 4.6 needs to be modified as follows.

LEMMA 4.6*. *For any $\delta>0$, there exists some $\varepsilon_0>0$ such that for any $k$,*

$$(81)\qquad \mathbb{P}_x(\tau_{\varepsilon,k}=N_k(x)\text{ and }\tau^{(k)}\leq\tau^{(j)}\text{ for all }j\neq k)\geq\varepsilon^{1/6+\delta}$$

*for all $\varepsilon\leq\varepsilon_0$, for all $x\in\Gamma$.*

Equation (81) implies, in particular, that

$$\mathbb{P}_x(\tau^{(k)}\leq\tau^{(j)}\text{ for all }j\neq k)\geq\varepsilon^{1/6+\delta}.$$

The argument of Section A.2 now gives

$$\sup_{x \in \overline{\gamma}} (\mathrm{Var}(P_1^n(x, dy) - \mu(dy))) \leq rc^{n\varepsilon^{1/6+\delta}},$$

(82)

$$\sup_{x \in \gamma} (\mathrm{Var}(P_2^n(x, dy) - \nu(dy))) \leq rc^{n\varepsilon^{1/6+\delta}}$$

for some $r > 0$ and $0 < c < 1$, which suffices for the purposes of Section 2.

To derive (81), we proceed as in the proof of Lemma 4.6. We split $N_k(x) = n_1 + n_2$ as in Section 5. Let us consider the more complicated case when $n_1 > 0$. The bad times are now defined by the condition

$$d(Y_n, B_j) \leq \sqrt{n}\varepsilon^{1/2-\alpha_1} \qquad \text{for some } j.$$

The contribution of the bad times is estimated as before, except that now we have "very bad times" when

$$d(X_n, B_j) \leq \varepsilon^{1/2-\delta} \qquad \text{for some } j \neq k$$

because then the orbit can be sucked into $U_j$. Let $n_3 = [\varepsilon^{-(1/2-2\delta)/(3/2+\delta)}]$. Observe that, with probability $1 - O(\varepsilon^R)$, the time difference between two very bad times corresponding to the same $j$ is at least $n_3$. Indeed, the probability that both $m_1$ and $m_2$ are very bad is $O(\varepsilon^R)$ unless $m = m_2 - m_1$ satisfies

$$\|m\rho\| < \sqrt{m}\varepsilon^{1/2-2\delta}.$$

This is impossible for $m \leq n_3$ since $m^{-(1+\delta)} \leq \|m\|$, due to Lemma 5.1. We can estimate $\mathbb{P}(\tau^{(j)} > n_3 \text{ for all } j \neq k)$ from below by a constant since there is at most one very bad time for each $j$ before $n_3$ and even if the orbit passes near a saddle point, it avoids $U_j$ with positive probability. Using arguments similar to those in the proof of Lemma 5.4, the conditional probability of having very bad time between $n_3$ and $n_1$, given that $\tau^{(j)} > n_3$ for all $j \neq k$, can be bounded by

$$\frac{n_1}{n_3} \frac{\varepsilon^{1/2-\delta}}{\sqrt{n_3\varepsilon}} = o\left(\frac{\sqrt{n_2}}{\sqrt{n_1}}\right).$$

Thus, we have the following analogue of Part (a) of Lemma 5.4.

(a*) There is a positive constant $k_1$ such that for all sufficiently small $\varepsilon$ and all $x$ such that $n_1 \neq 0$, we have

$$\mathbb{P}_x(\|\theta(X_{n_1}) + \rho n_2\| \leq \sqrt{n_2\varepsilon} \text{ and } \tau^{(j)} > n_1 \text{ for all } j \neq k) \geq \frac{k_1\sqrt{n_2}}{\sqrt{n_1}}.$$

Finally, $\mathbb{P}(\tau^{(j)} > N_k(x) \text{ for all } j \neq k | \tau^{(j)} \geq n_1 \text{ for all } j \neq k)$ can be bounded from below by a constant similarly to $\mathbb{P}(\tau^{\widehat{(j)}} > n_3 \text{ for all } j \neq k)$ and thus an analogue of Part (b) of Lemma 5.4 also remains valid.

The rest of the proof of (81) proceeds as in Section 5.

Finally, the result remains true if some of the periodic components contain more than one critical point. In this case, the edges $I_k$ should be replaced by subgraphs $\mathbb{G}_k$ which contain one common vertex, corresponding to the entire ergodic component. The other vertices of $\mathbb{G}_k$ correspond to the critical points of $H$ inside $U_k$ and the gluing conditions on those vertices are given in [7].

**Acknowledgments.** We are grateful to Prof. M. Freidlin for introducing us to this problem and for many useful discussions.

DEPARTMENT OF MATHEMATICS
UNIVERSITY OF MARYLAND
COLLEGE PARK, MARYLAND 20742
USA
E-MAIL: dmitry@math.umd.edu
        koralov@math.umd.edu